\documentclass[12pt,a4paper]{amsart}
\usepackage{geometry}
\usepackage{latexsym}
\usepackage{bm}
\usepackage{amsmath}
\usepackage{amssymb}
\usepackage{amsmath,amscd}
\usepackage{tabls}
\usepackage{url}
\usepackage[utf8]{inputenc}
\usepackage{amsfonts}
\usepackage{amssymb}
\usepackage{mathrsfs}
\usepackage{amsmath}
\usepackage{amsthm}
\usepackage{amscd}
\usepackage{fancyhdr}
\usepackage{enumerate}
\usepackage{paralist}
\usepackage{xypic}
\usepackage{graphicx}
\usepackage{cite}
\usepackage{lipsum}
\usepackage{footmisc}
\usepackage[all,cmtip]{xy}

\theoremstyle{plain}

\newtheorem{Def}[equation]{Definition}
\newtheorem{Thm}[equation]{Theorem}
\newtheorem{lem}[equation]{Lemma}
\newtheorem{prop}[equation]{Proposition}
\newtheorem{rem}[equation]{Remark}
\newtheorem{exa}[equation]{Example}
\newtheorem{conj}[equation]{Conjecture}

\numberwithin{equation}{section}

\title{Multiple Dedekind Symbols}

\author{Zhongyu Jin}
\address{Z. Jin \\School of Mathematical Sciences,
        Peking University,
         Beijing, China.}
\email{zyjin@pku.edu.cn}

\author{Qingchun Tian}
\address{Q. Tian \\School of Mathematical Sciences,
        Peking University,
         Beijing, China.}
\email{yhtian@math.pku.edu.cn}

\author{Yuan Zhou}
\address{Y. Zhou \\School of Mathematical Sciences,
        Peking University,
         Beijing, China.}
\email{zy93@pku.edu.cn}

\begin{document}
\maketitle

\begin{abstract}
In this paper we study multiple Dedekind symbols and the associated multiple reciprocity functions. There is a bijection between the two sets of them after a normalization. By this bijection we define products of multiple reciprocity functions, and study the relationship to the shuffle property. We construct and calculate shuffled multiple Dedekind symbols and shuffled multiple reciprocity functions from modular forms by regularized iterated integrals. Also we give a decomposition for them.
\end{abstract}

\section{Introduction}

In this paper we study the multiple Dedekind symbols, which generalize the Dedekind symbols of Fukuhara \cite{S.F4} and help understand holomorphic multiple modular values of Brown \cite{F.B}.

\subsection{Backgrounds}

The classical Dedekind sum is defined by
$$s(p,q)=\sum_{u=1}^{q-1}uq^{-1}(upq^{-1}-[upq^{-1}]-\frac{1}{2}),$$
where $p\in\mathbb{Z}$, $q\geq1$, and $[x]$ is the greatest integer in $x\in \mathbb{Q}$. It reflects modular properties of the Dedekind eta function, one can find more details in \cite{R.K}. Apostol \cite{T.A} defined generalized Dedekind sum and set up their relationship to modular forms. Following Apostol, Fukuhara \cite{S.F4} defined (generalized) Dedekind symbols and reciprocity functions with values in a commutative group.

For a general commutative group $G$, a function $D:U\rightarrow G$ is called a $G$-valued Dedekind symbol if it satisfies
$$D(p,-q)=D(-p,q), D(p,q)=D(p,p+q),$$
where $U=\{(p,q)\in\mathbb{Z}\times\mathbb{Z};\text{gcd}(p,q)=1\}$. We call $F(p,q)=D(p,q)D^{-1}(-q,p)$ the associated function of $D(p,q)$ (here $D^{-1}(p,q)$ denotes the inverse of $D(p,q)$ in $G$). The classical Dedekind sum can be regarded as an odd Dedekind symbol:
$$s(p,-q)=-s(p,q)\ \text{and}\ s(p,q)=s(p,p+q).$$
Also one can define $G$-valued reciprocity functions $F:U\rightarrow G$ satisfying that
$$F(p,-q)=F(-p,q), F(p,q)+F(-q,p)=0, F(p,p+q)+F(p+q,q)=F(p,q).$$
Then the function associated to a Dedekind symbol is a reciprocity function. This actually gives a bijection between the set of $G$-valued (normalized) Dedekind symbols and the set of $G$-valued reciprocity functions.

In particular, when $G$ is the group $\mathbb{C}$, Fukuhara totally determined the $\mathbb{C}$-valued Dedekind symbols with Laurent reciprocity law (i.e. ones whose associated functions are in Laurent polynomial form of $p$ and $q$) in \cite{S.F1}, \cite{S.F2}, \cite{S.F4} and \cite{S.F5}. All such functions can be constructed from modular forms. For instance, given a cusp form $f$ of weight $k\geq4$, define
\begin{equation}\label{length 1 DS of ES}
D_{f}(p,q)=\int_{q/p}^{i\infty}f(\tau)(p\tau-q)^{k-2}d\tau.
\end{equation}
Then $D_{f}(p,q)$ is a Dedekind symbol with Laurent reciprocity law, and the associated reciprocity function $F_{f}(p,q)$ is the period polynomial of $f$ (up to a constant term) by regarding $p$ and $q$ as variables. We will recall the construction for general modular forms in Section $2$.

Thus we have isomorphisms among the space of $\mathbb{C}$-valued odd (resp. even) Dedekind symbol with Laurent reciprocity law, the space of $\mathbb{C}$-valued odd (resp. even) reciprocity functions in Laurent polynomial form, and the space of cusp forms (resp. modular forms). These isomorphisms give a restatement of the Eichler-Shimura isomorphism. By this observation we can explain periods of modular forms by Dedekind symbols. As an example, consider the Dedekind symbol of the Eisenstein series of weight $2k\geq4$, take $q=1$ and $p\geq 2$, denote by $\xi=e^{2\pi iq/p}$, then we have the following reciprocity law:
\begin{equation}\label{length 1 reciprocity law}
\begin{split}
& \sum_{n=1}^{p-1}nLi_{2k-1}(\xi^{n})\\
= &-\frac{(2\pi i)^{2k-1}}{2} \sum_{n=-1,\ odd}^{2k-1}\frac{B_{n+1}B_{2k-n-1}}{(n+1)!(2k-n-1)!}p^{1-n}-\frac{(2\pi i)^{2k-1}B_{2k}}{4k(2k-2)!}p^{-2k+2}\\
& +\frac{\zeta(2k-1)}{2}(p^{-2k+3}-p),
\end{split}
\end{equation}
where $Li_{k}$ is the polylogarithm, $B_{k}$ is the $k$-th Bernoulli number and $\zeta(2k-1)$ is the zeta value. One can find more details in \cite{S.F2}.

As the generalization of the periods of modular forms, Manin \cite{Y.M1} defined iterated integrals of cusp form, Brown showed that such iterated integrals can be regarded as periods of his expected mixed modular motives, and a more accessible generalization is Brown's (holomorphic) multiple modular values defined in \cite{F.B}. It is natural to ask how to study such periods by Dedekind symbols.

In fact, for non-necessarily commutative group $G$, Manin \cite{Y.M3} define $G$-valued Dedekind symbols and reciprocity functions in the above way, and there is also a one-to-one correspondence between (normalized) Dedekind symbols and reciprocity functions.

In particular, let $G$ be the non-commutative group $\mathbb{C}\langle\langle \mathcal{A}\rangle\rangle^{\times}$ with $\mathcal{A}$ a set of symbols associated to cusp forms, Manin \cite{Y.M3} constructed $G$-valued Dedekind symbols by the iterated integrals of these cusp forms.

Manin's construction is our prototype of multiple Dedekind symbols, basing on which we define abstract multiple Dedekind symbols. Also we construct $\mathbb{C}$-valued ones from general modular forms, we need regularized iterated integrals of modular forms defined by Brown \cite{F.B}, and a construction of Bruggeman and Choie \cite{R.B}. The construction in \cite{R.B} was used to determine the non-commutative cohomology $H^{1}(SL_{2}(\mathbb{Z}),\mathcal{N}(A))$, where $\mathcal{N}(A)$ is the unit group of the non-commutative ring of formal power series in non-commutative variables in the set $A$, coefficients as holomorphic functions on the lower complex half plane.

\subsection{Outlines}

This paper is organized as below. In section $2$, we recall main results of Fukuhara and Manin. Then in Section $3$ and Section $4$, we study abstract multiple Dedekind symbols and multiple reciprocity functions.

In Section $3$, we define multiple Dedekind symbols. For a commutative unitary ring $R$, an $R$-valued multiple Dedekind symbol based on an indexed set $\mathcal{A}$ is a Dedekind symbol
$$D_{\mathcal{A}}:U\rightarrow R\langle\langle A\rangle\rangle^{\times}; D_{\mathcal{A}}(p,q)=1+\sum_{\emptyset\neq w\in \mathcal{A}^{*}}D^{w}(p,q)w,$$
where $\mathcal{A}^{*}$ is the set of words of $\mathcal{A}$. For any word $w\in \mathcal{A}^{*}$, we call $D^{w}(p,q)$ the $w$-component of $D_{\mathcal{A}}(p,q)$, which review extra information except the general properties of Dedekind symbols. With the help of Fukuhara's and Manin's work, we give the following result and briefly recall its proof.

\begin{Thm}
There is a bijection between the set of normalized multiple Dedekind symbols and the set of multiple reciprocity functions.
\end{Thm}

Moreover, we show that when we delete the point $(p,q)\in U$ satisfying that $pq=0$ from $U$, we can still define similar functions, which are called almost Dedekind symbols and reciprocity functions in this paper, and there is still a one-to-one correspondence in this case.

In Section $4$, we study a kind of product between multiple reciprocity functions, this product is induced by the product of multiple Dedekind symbols. Also we consider the shuffle property and its relationship to the product.

\begin{Thm}\label{Intro: shuffle and product}
The following statements hold:

$(1)$. A normalized multiple Dedekind symbol is shuffled if and only if its associated multiple reciprocity function is shuffled.

$(2)$. The inverse and products of shuffled multiple Dedekind symbols, and furthermore multiple reciprocity functions, are still shuffled.
\end{Thm}

In section $5$, we recall and generalize the regularized iterated integrals of Brown \cite{F.B}, and construct $\mathbb{C}$-valued shuffled multiple Dedekind symbols and multiple reciprocity functions from modular forms. More precisely, denote by $\mathcal{M}_{k}(\Gamma)$ the space of modular forms for $\Gamma=SL_{2}(\mathbb{Z})$ of weight $k$. Given $\mathcal{A}=\{A_{1},A_{2},\cdots\}$ as the indexed set, associate every element $A_{i}\in \mathcal{A}$ to an even integer $w_{i}\geq2$.  For any word $B=A_{1}\cdots A_{l}$, denote by $w(B)=w_{1}+\cdots+w_{l}$. Define
$$\mathcal{B}(\mathcal{A})=\{B\in\mathcal{A}^{*};B\neq \emptyset\}\ \text{and}\ \mathcal{M}(\mathcal{A},\Gamma)=\prod_{B\in \mathcal{B}(\mathcal{A})}\mathcal{M}_{w(B)+2}(\Gamma).$$
For any element $\mathbf{h}\in \mathcal{M}(\mathcal{A},\Gamma)$, denote by $\mathbf{h}(B)$ the component of $\mathbf{h}$ at the factor $B\in \mathcal{B}(\mathcal{A})$, it is a modular form of weight $w(B)+2$ for $\Gamma$. Denote by $I_{\mathbf{h}}(\delta_{1},\delta_{2})(p,q)$ the regularized iterated integral of
$$\Omega_{\mathbf{h}}(\tau)=\sum_{B\in\mathcal{B}(\mathcal{A})}B\mathbf{h}(B)(p\tau-q)^{w(B)}d\tau$$
from the tangential base point $\delta_{1}$ to $\delta_{2}$. For more details see Definition $5.2$.

\begin{Thm}
Given $\mathbf{h}\in\mathcal{M}(\mathcal{A},\Gamma)$, the regularized iterated integral
$$D_{\mathbf{h}}(p,q)=I_{\mathbf{h}}(\overrightarrow{q/p}_{\infty},\overrightarrow{\infty}_{q/p})(p,q)$$
is a $C$-valued shuffled multiple Dedekind symbol, and
$$F_{\mathbf{h}}(p,q)=I_{\mathbf{h}}(\overrightarrow{q/p}_{\infty},\overrightarrow{q/p}_{0})(p,q)\ \text{and}\ E_{\mathbf{h}}(p,q)=\text{exp}(\sum_{A_{i}\in\mathcal{A}}A_{i}a_{\mathbf{h}(A_{i})}(0)/pq)$$
are $\mathbb{C}$-valued shuffled multiple reciprocity functions. The reciprocity property is
$$D_{\mathbf{h}}(p,q)E_{\mathbf{h}}(p,q)D_{\mathbf{h}}^{-1}(-q,p)=F_{\mathbf{h}}(p,q).$$
\end{Thm}

In the same section, we calculate length two components of the multiple Dedekind symbols of Eisenstein series. Also we give the following decomposition of $\mathbb{C}$-valued multiple reciprocity functions, which help us determine their real or imaginary parts (the product here is the one in Theorem \ref{Intro: shuffle and product}). Although we only consider multiple Dedekind symbols and multiple reciprocity functions of modular forms, we actually proved the following theorem:

\begin{Thm}
Any shuffled multiple reciprocity function can be decomposed into a (maybe infinite) product of commutative valued reciprocity functions, and the corresponding result holds for shuffled multiple Dedekind symbols.
\end{Thm}

We explicitly construct these commutative valued reciprocity functions in the proof, they are unique after giving an order to words of $\mathcal{A}$. As a consequence, one deduce the problem that determining real (or imagine) part of $F_{\mathbf{h}}(p,q)$ as above to that determining ones of some reciprocity functions.

In the final section, we show how to define $\mathbb{C}$-valued (multiple) Dedekind symbols and (multiple) reciprocity functions for congruence subgroups of $SL_{2}(\mathbb{Z})$ of finite index, and give the analogue results. Also we calculate the Dedekind symbol of Eisenstein series.

\section{Dedekind Symbols and Reciprocity Functions}

In this section we recall the work of Fukuhara \cite{S.F4} and Manin \cite{Y.M3} about Dedekind symbols and reciprocity functions. Denote by
$$U=\{(p,q)\in\mathbb{Z}\times\mathbb{Z};\text{gcd}(p,q)=1\}.$$
Let $G$ be a (non-necessarily commutative) group. As mentioned in the introduction:

\begin{Def}
A $G$-valued Dedekind symbol is a function $D:U\rightarrow G$ satisfying that
$$D(p,-q)=D(-p,q), D(p,q)=D(p,p+q).$$
The function $f(p,q)=D(p,q)D^{-1}(-q,p)$ is called the associated function of $D(p,q)$.
\end{Def}

\begin{Def}
A $G$-valued reciprocity function is a function $F:U\rightarrow G$ satisfying that
$$F(p,-q)=F(-p,q), F(p,q)F(-q,p)=1, F(p,p+q)F(p+q,q)=F(p,q).$$
\end{Def}

\begin{rem}
The definitions given by Fukuhara \cite{S.F4} and Manin \cite{Y.M3} are not the same but differ little, we would like to use Manin's ones for convenience.
\end{rem}

For a given reciprocity function $F(p,q)$, consider the function
$$H:U\rightarrow G; H(p,q)=F(p,-q).$$
If $F=H$, we call $F$ an even reciprocity function. If $G$ is commutative and $F=H^{-1}$, we call $F$ an odd reciprocity function (if $G$ is non-commutative, $H$ is not a reciprocity function in general). When $G$ is commutative, one can write a reciprocity function as the sum of its even and odd parts.

The function $F(p,q)$ associated to a $G$-valued Dedekind symbol is naturally a $G$-valued reciprocity function. Since $F(p,q)=D(p,q)D^{-1}(-q,p)$, we have
$$F(p,-q)=D(p,-q)D^{-1}(q,p)=D(-p,q)D^{-1}(-q,-p)=F(-p,q),$$
$$F(p,q)F(-q,p)=D(p,q)D^{-1}(-q,p)D(-q,p)D^{-1}(-p,-q)=1,$$
$$F(p,p+q)F(p+q,q)=D(p,p+q)D^{-1}(-p-q,p)D(p+q,q)D^{-1}(-q,p+q)=F(p,q).$$
Fukuhara \cite{S.F4} proved the following theorem for commutative group $G$ and Manin \cite{Y.M3} generalized it to the general case. We will briefly recall its proof and give a generalized result in Section $3$.

\begin{Thm}[Fukuhara \cite{S.F4} and Manin \cite{Y.M3}]\label{GDS bijective RF}
For any given $G$-valued reciprocity function, there exists a $G$-valued Dedekind symbol associated to it.
\end{Thm}

Notice that $\mathbb{C}\cup\{\infty\}$ is not a group under the natural addition. One may call $D(p,q)$ (resp. $F(p,q)$) a $\mathbb{C}$-valued Dedekind symbol (resp. reciprocity function) by regarding it as a formal series of $p$ and $q$ satisfying the desired properties. A more accurate way is provided in Section $3.2$, in where we show that one can delete the points $(p,q)$ in $U$ satisfying that $pq=0$ and do the similar definition. In this case we still have the desired properties.

\begin{rem}
We keep this assumption when we discuss $\mathbb{C}$-valued multiple Dedekind symbols and multiple reciprocity functions below.
\end{rem}

Denote by $\mathcal{M}_{2k}(\Gamma)$ the complex vector space of modular forms for $\Gamma=SL_{2}(\mathbb{Z})$ of weight $2k$. For any $f\in \mathcal{M}_{2k}(\Gamma)$, define
\begin{equation*}
\begin{split}
& D_{f}(p,q)=\int_{\tau_{0}}^{i\infty}(f(\tau)-a_{f}(0))(p\tau-q)^{w}d\tau+\int_{q/p}^{\tau_{0}}(f(\tau)-\frac{a_{f}(0)}{(p\tau-q)^{w+2}})(p\tau-q)^{w}d\tau\\
& \ \ \ \ \ \ \ \ \ \ \ \ \ \ \ -a_{f}(0)(\frac{1}{(w+1)p}(p\tau_{0}-q)^{w+1}+\frac{1}{p}(p\tau_{0}-q)^{-1}),
\end{split}
\end{equation*}
where $a_{f}(0)$ is the Fourier constant term of $f$, $\tau_{0}$ is an arbitrary point on the upper half complex plane and $w=2k-2$. The function $D_{f}(p,q)$ is well-defined as mentioned by Zagier \cite{D.Z} and it is a $\mathbb{C}$-valued Dedekind symbol by Fukuhara \cite{S.F4}. When $f$ is a cusp form, it's the function in Formula (\ref{length 1 DS of ES}).

\begin{rem}\label{extra homo term}
When $f(\tau)=E_{2k}(\tau)$ is the Eisenstein series of weight $2k$, the associated reciprocity function of $D_{f}(p,q)$ is a sum of a homogeneous Laurent polynomial of degree $2k-2$ and an extra term $c/pq$ with $c$ a constant. One may find why the extra term exists in the proof of Theorem \ref{modular forms to MRF}.
\end{rem}

Define maps
$$\phi_{2k}:\mathcal{M}_{2k}(\Gamma)\rightarrow \mathcal{D};f\mapsto D_{f},$$
$$\psi:\mathcal{D}\rightarrow\mathcal{F}; \psi(D)(p,q)=D(p,q)-D(-q,p),$$
where $\mathcal{D}$ is the space of $\mathbb{C}$-valued Dedekind symbols and $\mathcal{F}$ is the space of $\mathbb{C}$-valued reciprocity functions. The function $\phi_{2k}$ is well-defined, denote by $\phi_{2k}^{+}$ and $\phi_{2k}^{-}$ the even and odd part of $\phi_{2k}$ respectively as in Fukuhara \cite{S.F4}.

Denote by $\mathcal{FL}_{2k}^{+}$ the space of even reciprocity functions in homogeneous polynomial form of degree $2k-2$, and by $\mathcal{FL}_{2k,-1}^{-}$ the space of odd reciprocity functions in the form $G(p,q)/pq$, where $G(p,q)$ is a homogeneous polynomial of degree $2k$ (up to an extra term as in Remark \ref{extra homo term}). Fukuhara \cite{S.F4} proved the following theorem, which restate the Eichler-Shimura isomorphism.

\begin{Thm}[Fukuhara \cite{S.F4}]\label{DS, RF and MF}
The map $\psi$ is an isomorphism of $\mathbb{C}$-vector spaces. For any $k\geq 2$, the map $\phi_{2k}$ is injective. Furthermore, we have
$$\psi\circ\phi_{2k}^{+}(\mathcal{M}_{2k}(\Gamma))=\mathcal{FL}_{2k}^{+}\ \text{and}\ \psi\circ\phi_{2k}^{-}(\mathcal{M}_{2k}(\Gamma))=\mathcal{FL}_{2k,-1}^{-}.$$
\end{Thm}

\section{Multiple Dedekind Symbols and Multiple Reciprocity Functions}

Let $\mathcal{A}=\{A_{1},A_{2},\cdots\}$ be an indexed set (maybe infinite), we call $u=A_{i_{1}}\cdots A_{i_{r}}$ a word of $\mathcal{A}$ for $A_{i_{j}}\in \mathcal{A}$, and $l(w)=r$ the length of $w$. The empty word $\emptyset$ is also regarded as a word of $\mathcal{A}$, and its length is $0$. Denote by $\mathcal{A}^{*}$ the set of words of $\mathcal{A}$, there is a natural contraction product on $\mathcal{A}^{*}$, i.e., for $u=A_{1}\cdots A_{r}$ and $v=B_{1}\cdots B_{s}$,
$$uv=A_{1}\cdots A_{r}B_{1}\cdots B_{s}\in \mathcal{A}^{*}.$$

In the following part of this paper, we always assume that $R$ is a commutative unitary ring. Denote by $R\langle\langle \mathcal{A}\rangle\rangle$ the ring of formal power series of words of $A$ with coefficients in $R$ and by $R\langle\langle \mathcal{A}\rangle\rangle^{\times}$ the group of invertible elements in it. Then $R\langle\langle \mathcal{A}\rangle\rangle^{\times}$ is a group, and is non-commutative when $|\mathcal{A}|\geq 2$.

\subsection{Fundamental results}

Recall that $U=\{(p,q)\in\mathbb{Z}\times\mathbb{Z};\text{gcd}(p,q)=1\}$.

\begin{Def}
An $R$-valued multiple Dedekind symbol based on the set $\mathcal{A}$ is an $R\langle\langle \mathcal{A}\rangle\rangle^{\times}$-valued Dedekind symbol in the form
$$D_{\mathcal{A}}:U\rightarrow R\langle\langle \mathcal{A}\rangle\rangle^{\times}; D_{\mathcal{A}}(p,q)=1+\sum_{\emptyset\neq w\in \mathcal{A}^{*}}D^{w}(p,q)w.$$
It is invertible since the constant term of $D_{\mathcal{A}}(p,q)$ is $1$. We call the multiple Dedekind symbol $D_{\mathcal{A}}$ normalized if $D_{\mathcal{A}}(1,1)=1$.
\end{Def}

If we denote by
$$F_{\mathcal{A}}(p,q)=D_{\mathcal{A}}(p,q)D_{\mathcal{A}}^{-1}(-q,p)=1+\sum_{\emptyset\neq w\in \mathcal{A}^{*}}F^{w}(p,q)w$$
the associated function, then the definition is equivalent to the following conditions: for any non-empty word $w\in \mathcal{A}^{*}$,

$(MDS\ 1)$. $D^{w}(p,-q)=D^{w}(-p,q)$.

$(MDS\ 2)$. $D^{w}(p,q)=D^{w}(p,p+q)$.\\
Also, we have

$(MDS\ 3)$. $D^{w}(p,q)=\sum\limits_{uv=w}F^{u}(p,q)D^{v}(-q,p),$\\
where $D^{w}(p,q)=1$ when $w=\emptyset$. We call $D^{w}(p,q)$ the $w$-component of $D_{\mathcal{A}}(p,q)$ and $l(w)$ its length. $D_{\mathcal{A}}(p,q)$ is normalized if and only if $D^{w}(1,1)=0$ for all $w\neq\emptyset$.

\begin{lem}\label{unique nMDS}
For any $R$-valued multiple Dedekind symbol $D_{\mathcal{A}}$ based on $\mathcal{A}$ associated to the function $F_{\mathcal{A}}$, there is a unique normalized multiple Dedekind symbol $\widetilde{D}_{\mathcal{A}}$ whose associated function is still $F_{\mathcal{A}}$.
\end{lem}

\noindent{\bf Proof:}
Suppose that $D_{\mathcal{A}}(p,q)=1+\sum\limits_{\emptyset\neq w\in \mathcal{A}^{*}}D^{w}(p,q)w$, directly from the definition, the constant function
$$D_{\mathcal{A}}^{c}(p,q)=D_{\mathcal{A}}(1,1)=1+\sum\limits_{\emptyset\neq w\in \mathcal{A}^{*}}D^{w}(1,1)w$$
is also a multiple Dedekind symbol, and its associated function is the constant function $1$. We can define the function $\widetilde{D}_{\mathcal{A}}:U\rightarrow R\langle\langle \mathcal{A}\rangle\rangle^{\times}$ by
$$\widetilde{D}_{\mathcal{A}}(p,q)=D_{\mathcal{A}}(p,q)D_{\mathcal{A}}^{c}(p,q)^{-1}.$$
Obviously the function $\widetilde{D}_{\mathcal{A}}$ is an $R$-valued normalized multiple Dedekind symbol based on $\mathcal{A}$, whose associated function is $F_{\mathcal{A}}$.

Next, we show the uniqueness. Suppose that $\widetilde{D}_{\mathcal{A}}$ and $\widetilde{D}_{\mathcal{A}}'$ are two normalized multiple Dedekind symbols associated to the function $F_{\mathcal{A}}$. Denote by
$$K_{\mathcal{A}}(p,q)=\widetilde{D}_{\mathcal{A}}(p,q)^{-1}\widetilde{D}_{\mathcal{A}}'(p,q).$$
It is obvious to see
$$K_{\mathcal{A}}(p,q)=K_{\mathcal{A}}(p,p+q),$$
and
\begin{equation*}
\begin{split}
K_{\mathcal{A}}(p,q)&=\widetilde{D}_{\mathcal{A}}(p,q)^{-1}\widetilde{D}'_{\mathcal{A}}(p,q)\\
&=\widetilde{D}_{\mathcal{A}}(-q,p)^{-1}F_{\mathcal{A}}(p,q)^{-1}F_{\mathcal{A}}(p,q)\widetilde{D}'_{\mathcal{A}}(-q,p)\\ &=K_{\mathcal{A}}(-q,p).
\end{split}
\end{equation*}
As a consequence, $K_{\mathcal{A}}$ is a constant function. Since the multiple Dedekind symbols $\widetilde{D}_{\mathcal{A}}$ and $\widetilde{D}'_{\mathcal{A}}$ are both normalized, we have $K_{\mathcal{A}}(p,q)=K_{\mathcal{A}}(1,1)=1$, and thus $\widetilde{D}_{\mathcal{A}}=\widetilde{D}'_{\mathcal{A}}$.
$\hfill\Box$\\

\begin{Def}
An $R$-valued multiple reciprocity function based on the set $\mathcal{A}$ is an $R\langle\langle \mathcal{A}\rangle\rangle^{\times}$-valued reciprocity function in the form
$$F_{\mathcal{A}}:U\rightarrow R\langle\langle \mathcal{A}\rangle\rangle^{\times}; F_{\mathcal{A}}(p,q)=1+\sum_{\emptyset\neq w\in \mathcal{A}^{*}}F^{w}(p,q)w.$$
\end{Def}

Similarly, the definition of $R$-valued multiple reciprocity functions based on $\mathcal{A}$ is equivalent to the following conditions: for any non-empty word $w\in \mathcal{A}^{*}$,

$(MRF\ 1)$. $F^{w}(p,-q)=F^{w}(-p,q)$.

$(MRF\ 2)$. $\sum\limits_{uv=w}F^{u}(p,q)F^{v}(-q,p)=0, w\neq\emptyset,$

$(MRF\ 3)$. $\sum\limits_{uv=w}F^{u}(p,p+q)F^{v}(p+q,q)=F^{w}(p,q)$,\\
where $F^{w}(p,q)=1$ when $w=\emptyset$, and we also call $F^{w}(p,q)$ the $w$-component of $F_{A}(p,q)$ and $l(w)$ its length. For any $w\in\mathcal{A}^{*}$ of length $1$, $F^{w}(p,q)$ is an $R$-valued reciprocity function.

\begin{rem}
For word $w\in \mathcal{A}^{*}$ with $l(w)\geq2$, the component $F^{w}$ is not an $R$-valued reciprocity function in general.
\end{rem}

When there is no confusion of the choice of the ring $R$ and the set $\mathcal{A}$, we would like to call $D_{\mathcal{A}}$ and $F_{\mathcal{A}}$ the multiple Dedekind symbol and the multiple reciprocity function respectively for simplicity. In section $4$, we show that one can regard a general Dedekind symbol (resp. reciprocity function) as a multiple one. Combining with Theorem \ref{GDS bijective RF} and Lemma \ref{unique nMDS}, we have

\begin{Thm}\label{MDS bijective MRF}
There is a bijection between the set of normalized multiple Dedekind symbols and the set of  multiple reciprocity functions.
\end{Thm}

We briefly recall how to construct a normalized multiple Dedekind symbols from a given multiple reciprocity function. One can find more detailed proof in \cite{S.F4} and \cite{Y.M3}. In the next subsection, we give a generalized result.

Assume that $X_{0}, X_{1},\cdots$ is a set of (non-commutative) symbols, define polynomials $P^{(n)},Q^{(n)}\in\mathbb{Z}[X_{0},\cdots,X_{n}]$ inductively to be:
$$P^{(n+1)}(X_{0},\cdots,X_{n+1})=Q^{(n)}(X_{1},\cdots,X_{n+1}),$$
$$Q^{(n+1)}(X_{0},\cdots,X_{n+1})=X_{0}Q^{(n)}(X_{1},\cdots,X_{n+1})-P^{(n)}(X_{1},\cdots,X_{n+1}),$$
$$P^{(0)}(X_{0})=1, Q^{(0)}(X_{0})=X_{0}.$$
Denote by $\langle X_{0},\cdots,X_{n}\rangle=Q^{(n)}(X_{0},\cdots,X_{n})/P^{(n)}(X_{0},\cdots,X_{n})$, then
$$\langle X_{0},\cdots,X_{n}\rangle=X_{0}-\frac{1}{X_{1}-\frac{1}{X_{2}-\cdots}},$$
where the lowest term in the continued fraction representation is $X_{n-1}-\frac{1}{X_{n}}$. By the definition of $P^{(n)}$ and $Q^{(n)}$, for any $n\in\mathbb{Z}_{\geq0}$ we have $\text{gcd}(P^{(n)},Q^{(n)})=1$. Given a multiple reciprocity function $F_{\mathcal{A}}(p,q)$, we construct the multiple Dedekind symbol $D_{\mathcal{A}}(p,q)$ as below.

Given $(p,q)\in U$, take a representation of $q/p=\langle a_{0},\cdots,a_{n}\rangle$ for some $n\in\mathbb{Z}_{\geq0}$ and $a_{i}\in\mathbb{Z}$. Denote by
$$q_{i}=Q^{(n-i)}(a_{i},\cdots,a_{n}),\ p_{i}=P^{(n-i)}(a_{i},\cdots,a_{n}),$$
then we have $q_{i}/p_{i}=\langle a_{i},\cdots,a_{n}\rangle$ and $\text{gcd}(p_{i},q_{i})=1$. Define
\begin{equation}\label{Construction of DS from RF}
D_{\mathcal{A}}(p,q)=\prod_{i=1}^{n}F_{\mathcal{A}}^{-1}(p_{i},q_{i}),
\end{equation}
which does not depend on the choice of the continued fraction representation of $q/p$. Then one finds that $D_{\mathcal{A}}$ is a normalized multiple Dedekind symbol associated to the given multiple reciprocity function $F_{\mathcal{A}}$.

\subsection{A generalization}

The independence of $D_{\mathcal{A}}(p,q)$ in equation (\ref{Construction of DS from RF}) on the representation of $q/p=\langle a_{0},\cdots,a_{n}\rangle$ follows from the following fact. Given a finite sequence of integers $\mathbf{a}=(a_{0},\cdots,a_{n})$, define three types of moves on it by

$(T1)$. Given $\varepsilon\in\{1,-1\}$ and $i\in\{0,\cdots,n-1\}$,
$$\mathbf{a}\mapsto \mathbf{b}_{1}=(a_{0},\cdots,a_{i}+\varepsilon,\varepsilon,a_{i+1}+\varepsilon,a_{i+2},\cdots).$$

$(T2)$. Given $i\in\{0,\cdots,n\}$ and a representation of $a_{i}=b+c$ for $b,c\in\mathbb{Z}$,
$$\mathbf{a}\mapsto \mathbf{b}_{2}=(a_{0},\cdots,a_{i-1},b,0,c,a_{i+1},\cdots,a_{n}).$$

$(T3)$. Given $\varepsilon\in\{1,-1\}$,
$$\mathbf{a}\mapsto \mathbf{b}_{3}=(a_{0},\cdots,a_{n-1},a_{n}+\varepsilon,\varepsilon).$$
One can find the proof of the following lemma in \cite{R.K}.

\begin{lem}\label{cont frac trans}
If two sequences $\mathbf{a}$ and $\mathbf{a}'$ differ by one of the above moves, then their respective continued fractions define the same rational number. Conversely, if we write a rational number in two different ways as a continued fraction, then the respective sequences differ by a finite sequence of moves above and their inverses.
\end{lem}

Also one can find the proof of the following theorem in \cite{S.F4} and \cite{Y.M3}.

\begin{Thm}\label{Indep of cont frac}
The function $D_{\mathcal{A}}(p,q)$ defined in Formula (\ref{Construction of DS from RF}) does not change when we do the moves $(T1)$, $(T2)$, $(T3)$ and their inverses. Thus $D_{\mathcal{A}}(p,q)$ depends only on $(p,q)$.
\end{Thm}

Theorem \ref{Indep of cont frac} holds for general Dedekind symbols. As stated in the Section $2$, we delete points $(p,q)$ of $U$ satisfying $pq=0$ from $U$, and consider
$$U'=\{(p,q)\in\mathbb{Z}\times\mathbb{Z}|\text{gcd}(p,q)=1,\ pq\neq 0\}.$$
Suppose that $G$ is a multiplicative group, non-necessarily commutative.

\begin{Def}
A G-valued almost Dedekind symbol is a function $D:U'\to G$ satisfying
$$D(p,q)=D(p,p+q),\ (p,q,p+q\neq0)\ \text{and}\ D(p,-q)=D(-p,q),\ (p,q\neq0).$$
We call $D$ normalized if $D(1,1)=1$.
\end{Def}

\begin{Def}
A G-valued almost reciprocity function is a function $F:U'\to G$ satisfying
$$F(p,q)=F(p,p+q)F(p+q,q),\ (p,q,p+q\neq0),$$
$$F(p,q)F(-q,p)=1,F(p,-q)=F(-p,q),\ (p,q\neq 0).$$
\end{Def}

One finds that a Dedekind symbol (resp. reciprocity function) restricted on $U'$ is an almost one. The $\mathbb{C}$-valued Dedekind symbols and reciprocity functions of modular forms (especially Eisenstein series) defined in Section $2$ are almost.

One can also do the normalization for almost Dedekind symbols as in Lemma \ref{unique nMDS}, and we have the following result:

\begin{Thm}\label{Almost DS bij RF}
Given a group $G$ as above, there is a one-to-one correspondence between the normalized $G$-valued almost Dedekind symbols and the $G$-valued almost reciprocity functions.
\end{Thm}

\noindent{\bf Proof:}
The difficulty is that one cannot use Fukuhara's and Manin's construction directly since we cannot get $F(1,1)=1$ from the definition. In order to solve this problem, we use a specialized continued fraction representation of $q/p$ as
\begin{equation}\label{Unique repre}
q/p=\langle a_{0},a_{1},\cdots,a_{n}\rangle,\ a_{i}\geq 2\ \text{for}\ 1\leq i\leq n.
\end{equation}
This representation is unique, and we will avoid using the move $(T3)$ below since from the proof of Theorem \ref{Indep of cont frac}, $(T1)$ and $(T2)$ works well in our case, but $(T3)$ does not since $F(1,1)$ may not equal to $1$.

By abusing the notations, still denote by $\mathcal{D}$ the set of normalized almost Dedekind symbols, by $\mathcal{F}$ the set of almost reciprocity function, and by $\psi:\mathcal{D}\to\mathcal{F}$ the map
$$\psi(D)(p,q)=D(p,q)D(q,-p)^{-1}.$$
It is obvious that $\psi(D)\in\mathcal{F}$.

Let $q/p=\left\langle a_0,a_1,\cdots,a_n\right\rangle$ be the continued fraction expansion, denote $(p,q)$ by $\left\langle a_0,a_1,\cdots,a_n\right\rangle$, and $(p_{i},q_{i})$ by $\langle a_{i},\cdots,a_{n}\rangle$, for convenience. We define the converse map $\delta:\mathcal{F}\to\mathcal{D}$ below.

If $q/p$ is a positive integer, we define $\delta(F)(p,q)=1$, and if $q/p$ is a negative integer, we define $\delta(F)(p,q)=F(1,1)^{-1}$. Assume that $q/p\in\mathbb{Q}\backslash\mathbb{Z}$, then it can be represented by $\left\langle a_0,a_1,\cdots,a_n\right\rangle$ as in formula (\ref{Unique repre}) with $a_i\geq 2$ for $1\leq i\leq n$, and we define
$$\delta(F)(p,q)=F(\left\langle a_1,\cdots,a_n\right\rangle)^{-1}F(\left\langle a_2,\cdots,a_n\right\rangle)^{-1}\cdots F(\left\langle a_n\right\rangle)^{-1}.$$
It is obvious that $\delta(F)\in\mathcal{D}$, we only need to prove $\delta\circ\psi=\text{id}$ and $\psi\circ\delta=\text{id}$.

First, we prove $\delta\circ\psi=\text{id}$. When $q/p$ in an integer, it is trivial, and we only consider the case that $q/p=\langle a_{0},\cdots,a_{n}\rangle\in\mathbb{Q}\backslash\mathbb{Z}$ as in Formula (\ref{Unique repre}). We have
\begin{equation*}
\begin{split}
& \delta\circ\psi(D)(\langle a_0,\cdots,a_n\rangle)\\
=& \psi(D)(\left\langle a_1,\cdots,a_n\right\rangle)^{-1}\cdots \psi(D)(\left\langle a_n\right\rangle)^{-1}\\
=& D(\left\langle0,a_1,\cdots,a_n\right\rangle)D(\left\langle a_1,a_2,\cdots,a_n\right\rangle)^{-1}\cdots D(\left\langle0,a_n\right\rangle)D(\left\langle a_n\right\rangle)^{-1}\\
=& D(\left\langle a_0,a_1,\cdots,a_n\right\rangle).
\end{split}
\end{equation*}
Thus one side holds.

Next we show $\psi\circ\delta=\text{id}$. If $q/p=n$ is a positive integer, it is trivial. If $q/p=-n$ is a negative integer, the following equation holds:
\begin{equation*}
\begin{split}
\psi\circ\delta(F)(-1,n)&=\delta(F)(-1,n)\delta(F)(n,1)^{-1}\\
&=F(1,1)^{-1}\delta(F)(\langle1,2,\cdots,2\rangle)^{-1}\\
&=F(-1,1)F(1,2)F(2,3)\cdots F(n-1,n)\\
&=F(-1,2)F(2,3)\cdots F(n-1,n)\\
&\qquad\vdots\\
&=F(-1,n)
\end{split}
\end{equation*}
since $F(\langle 2,\cdots,2\rangle)=F(k,k+1)$, where $2$ appears $k$ times in the left hand side.

If $q/p\in\mathbb{Q}\backslash\mathbb{Z}$, we have four conditions and we discuss them respectively. First, if $q/p=\left\langle a_0,\cdots,a_n\right\rangle>1$ and $a_i\geq 2$ for positive $i$, then we must have $a_0\geq2$ and
\begin{equation*}
\begin{split}
& \psi\circ\delta(F)(\left\langle a_0,\cdots,a_n\right\rangle)\\
=& \delta(F)(\left\langle a_0,\cdots,a_n\right\rangle)\delta(F)(\left\langle 0, a_0,\cdots,a_n\right\rangle)^{-1}\\
=& F(\left\langle a_1,\cdots,a_n\right\rangle)^{-1}\cdots F(\langle a_{n}\rangle)^{-1}F(\langle a_{n}\rangle)\cdots F(\left\langle a_1,\cdots,a_n\right\rangle)F(\left\langle a_0,\cdots,a_n\right\rangle)\\
=& F(\left\langle a_0,\cdots,a_n\right\rangle).
\end{split}
\end{equation*}

Then if $0<q/p=\left\langle a_0,\cdots,a_n\right\rangle<1$ with $a_i\geq 2$ for positive $i$, then we must have $a_0=1$. Notice that there exists at least one $i$ such that $a_{i}>2$, suppose that
$$q/p=\langle1,2,\cdots,2,a_k,\cdots,a_n\rangle,\ a_k\geq 3,$$
then we have
\begin{equation*}
\begin{split}
\psi\circ\delta(F)(\left\langle 1,2,\cdots,a_n\right\rangle)&=\delta(F)(\left\langle 1,2,\cdots,a_n\right\rangle)\delta(F)(\left\langle 0, 1,2,\cdots,a_n\right\rangle)^{-1}.
\end{split}
\end{equation*}
One only needs to prove
\begin{equation}\label{Case 2}
\delta(F)(\langle 0,1,2,\cdots,a_{n}\rangle)=F^{-1}(\langle1,2,\cdots,a_{n}\rangle)\cdots F^{-1}(\langle a_{n}\rangle).
\end{equation}
By Lemma \ref{cont frac trans} and Theorem \ref{Indep of cont frac}, if one can transform the specialized continued fraction representation $q/p=\langle a_{0}',\cdots,a_{m}'\rangle$ as in Formula (\ref{Unique repre}) to the continued fraction $\langle 0,1,2,\cdots,a_{n}\rangle$ by moves $(T1)$ and $(T2)$, then Equation (\ref{Case 2}) holds by the definition of $\delta$. One can do this by using the move $(T1)$ repeatedly as below:
\begin{equation*}
\begin{split}
\delta(F)(\left\langle 0,1,2,\cdots,a_n\right\rangle)&=\delta(F)(\left\langle 0,a_1-1,a_2\cdots,a_n\right\rangle)\\
&=\delta(F)(\left\langle 0,1,a_2\cdots,a_n\right\rangle)\\
&=\delta(F)(\left\langle -1,a_2-1,a_3,\cdots,a_n\right\rangle)\\
&\qquad\vdots\\
&=\delta(F)(\left\langle 1-k,a_k-1,\cdots,a_n\right\rangle).
\end{split}
\end{equation*}
By assumption, $a_{k}\geq3$, then the representation $\langle 1-k,a_{k}-1,\cdots,a_{n}\rangle$ is the one as formula (\ref{Unique repre}). Thus $\psi\circ \delta=\text{id}$ in this case.

If $-1<q/p=\left\langle a_0,\cdots,a_n\right\rangle<0$ with $a_i\geq 2$ for positive $i$, then we must have $a_0=0$ and
\begin{equation*}
\begin{split}
\psi\circ\delta(F)(\left\langle 0,a_1,\cdots,a_n\right\rangle)&=\delta(F)(\left\langle 0,a_1,\cdots,a_n\right\rangle)\delta(F)(\left\langle 0, 0,a_1,\cdots,a_n\right\rangle)^{-1}\\
& =\delta(F)(\left\langle 0,a_1,\cdots,a_n\right\rangle)\delta(F)(\left\langle a_1,\cdots,a_n\right\rangle)^{-1}\\
& =F(\langle 0,a_{1},\cdots,a_{n}\rangle).
\end{split}
\end{equation*}

Finally, if $q/p=\left\langle a_0,\cdots,a_n\right\rangle<-1$ with $a_i\geq 2$ for positive $i$, then we must have $a_0=-n\leq -1$ and
\begin{equation*}
\begin{split}
\psi\circ\delta(F)(\left\langle -n,a_1,\cdots,a_n\right\rangle)&=\delta(F)(\left\langle -n,a_1,\cdots,a_n\right\rangle)\delta(F)(\left\langle 0, -n,a_1,\cdots,a_n\right\rangle)^{-1}.
\end{split}
\end{equation*}
Similar to the argument in the second case, one can use the move $(T1)$ and $(T2)$ repeatedly and get
\begin{equation*}
\begin{split}
\delta(F)(\left\langle 0,-n,a_1,\cdots,a_n\right\rangle)&=\delta(F)(\left\langle 1,1,-n+1,a_1\cdots,a_n\right\rangle)\\
&=\delta(F)(\left\langle 1,2,1,-n+2,a_1\cdots,a_n\right\rangle)\\
&\qquad\vdots\\
&=\delta(F)(\left\langle 1,2,\cdots,2,1,0,a_1,\cdots,a_n\right\rangle)\\
&=\delta(F)(\left\langle 1,2,\cdots,2,1+a_1,\cdots,a_n\right\rangle).
\end{split}
\end{equation*}
Thus we have $\psi\circ\delta=\text{id}$ here. As a summary, the theorem holds.
$\hfill\Box$\\

For our later simplification of notations, we call all functions we discuss in the rest part of this paper, whether almost or not, (multiple) Dedekind symbols and (multiple) reciprocity functions, and write them as functions defined on $U$.

\section{Product and Shuffle Property}

In this section, we consider two further structures on multiple Dedekind symbols and multiple reciprocity functions. All the discussion and results in this section hold for almost ones.

\subsection{Products of multiple reciprocity functions}

Let $\mathcal{A}$ be an indexed set. In general, the natural product induced by the one in $R\langle\langle \mathcal{A}\rangle\rangle^{\times}$ of two $R$-valued multiple reciprocity functions $F_{\mathcal{A}}$ and $G_{\mathcal{A}}$ based on $\mathcal{A}$ is not a multiple reciprocity function based on $\mathcal{A}$.

However, since we have a one-to-one correspondence between normalized multiple Dedekind symbols and multiple reciprocity functions, and the product of two normalized multiple Dedekind symbols is still a normalized multiple Dedekind symbol, we can make the following definition.

Denote by $D_{\mathcal{A}}$ and $E_{\mathcal{A}}$ the associated normalized multiple Dedekind symbols of $F_{\mathcal{A}}$ and $G_{\mathcal{A}}$ respectively.

\begin{Def}
Define the product of $F_{\mathcal{A}}$ and $G_{\mathcal{A}}$ to be the multiple reciprocity function $F_{\mathcal{A}}\bullet G_{\mathcal{A}}$, whose associated normalized multiple Dedekind symbol is $D_{\mathcal{A}}E_{\mathcal{A}}$. Precisely for any $(p,q)\in U$,
\begin{equation*}
\begin{split}
(F_{\mathcal{A}}\bullet G_{\mathcal{A}})(p,q)& =D_{\mathcal{A}}(p,q)E_{\mathcal{A}}(p,q)E_{\mathcal{A}}^{-1}(-q,p)D_{\mathcal{A}}^{-1}(-q,p)\\
& = D_{\mathcal{A}}(p,q)G_{\mathcal{A}}(p,q)D_{\mathcal{A}}^{-1}(-q,p).
\end{split}
\end{equation*}
\end{Def}

When $F_{\mathcal{A}}$ is multiple reciprocity functions and $\mathcal{A}\subset \mathcal{A}'$, one can regard $F_{\mathcal{A}}$ as a multiple reciprocity function based on $\mathcal{A}'$ by letting $F^{w}=0$ for $w\notin \mathcal{A}^{*}$ (the corresponding result holds for multiple Dedekind symbols). Thus one can give the definition above when the indexed sets of the two multiple reciprocity functions are different.

In this and the following sections, all products of multiple reciprocity functions are the ones under $\bullet$ if there is no confusion. This product satisfies the following fundamental properties:

\begin{prop}
The product for multiple reciprocity functions is associated. The constant function $1$ is the unit, and any multiple reciprocity function $F_{\mathcal{A}}$ has a unique inverse under the product.
\end{prop}

\begin{rem}
The product is non-commutative in general.
\end{rem}

Given a multiple reciprocity function and the associated normalized multiple Dedekind symbol
$$F_{\mathcal{A}}=1+\sum_{\emptyset\neq w\in \mathcal{A}^{*}}F^{w}w,\ D_{\mathcal{A}}=1+\sum_{\emptyset\neq w\in \mathcal{A}^{*}}D^{w}w,$$
and another multiple reciprocity function
$$G_{\mathcal{A}}=1+\sum_{w\in \mathcal{A}^{*}}G^{w}w.$$
The $w$-components of $(F_{\mathcal{A}}\bullet G_{\mathcal{A}})$ is $F^{w}+G^{w}$ for any $w$ of length $1$. Also for any word $w=uv$ with $u,v\in \mathcal{A}$, the $w$-component of $(F_{\mathcal{A}}\bullet G_{\mathcal{A}})$ takes value
$$F^{uv}(p,q)+G^{uv}(p,q)+D^{u}(p,q)G^{v}(p,q)-G^{u}(p,q)D^{v}(-q,p)$$
at $(p,q)\in U$.

Moreover, according to the law $(MDS\ 3)$ and the discussion in Theorem \ref{Shuffled equivalence} below, for any positive integer $n$, the condition that $F^{w}=0$ for all $0<l(w)<n$ is equivalent to the condition that $D^{w}=0$ for all $0<l(w)<n$ (notice that $D_{\mathcal{A}}$ is normalized). As a consequence, we have the following result:

\begin{prop}
Denote by $m$ (resp. $n$) the biggest positive integer such that the length $k$ components of $F_{\mathcal{A}}$ (resp. $G_{\mathcal{A}}$) are all zero for $0<k<m$ (resp. $0<k<n$). Then the length $0<k<min\{m,n\}$ components of $(F_{\mathcal{A}}\bullet G_{\mathcal{A}})$ are all zero, and for any word $w$ of length $min\{m,n\}$, the $w$-component of $(F_{\mathcal{A}}\bullet G_{\mathcal{A}})$ is $F^{w}+G^{w}$.
\end{prop}

Finally, we give an example explaining that a Dedekind symbol, or a reciprocity function, can be naturally regarded as a multiple one. Suppose that $\mathcal{A}=\{a\}$ and that we are given an $R$-valued reciprocity function $F$. Define
$$F_{\mathcal{A}}:U\rightarrow R\langle\langle \mathcal{A}\rangle\rangle^{\times};F_{\mathcal{A}}(p,q)=e^{F(p,q)a}=1+\sum_{n\geq1}\frac{F(p,q)^{n}}{n!}a^{n},$$
where $a^{n}$ means the word $a\cdots a$ consisting of n $a$. Since there is only one element in the set $\mathcal{A}$, the product of words of $\mathcal{A}^{*}$ is commutative and obviously $F_{\mathcal{A}}$ is an $R$-valued multiple reciprocity function, and an $R\langle\langle \mathcal{A}\rangle\rangle^{\times}$-valued reciprocity function. Denote by $D_{\mathcal{A}}$ the associated normalized Dedekind symbol of $F_{\mathcal{A}}$, then
$$D_{\mathcal{A}}:U\rightarrow R\langle\langle \mathcal{A}\rangle\rangle^{\times};D_{\mathcal{A}}(p,q)=e^{D(p,q)a}=1+\sum_{n\geq1}\frac{D(p,q)^{n}}{n!}a^{n}$$
is an $R$-valued normalized multiple Dedekind symbol, and an $R\langle\langle \mathcal{A}\rangle\rangle^{\times}$-valued Dedekind symbol. The associated multiple reciprocity function of $D_{\mathcal{A}}$ is just the function $F_{\mathcal{A}}$. The function
$$G_{\mathcal{A}}(p,q)=e^{-F(p,q)a}=1+\sum_{n\geq 1}\frac{(-1)^{n}F(p,q)^{n}}{n!}q^{n}$$
is the inverse of $F_{\mathcal{A}}$ under both the product $\bullet$ and the product induced by the one on $R\langle\langle\mathcal{A}\rangle\rangle^{\times}$.

Conversely, if $\mathcal{A}=\{a\}$, for a normalized multiple Dedekind symbol $D_{\mathcal{A}}(p,q)$, we can define
$$D(p,q)=\log(D_{\mathcal{A}}(p,q))=D_{\mathcal{A}}(p,q)-1+\frac{1}{2}(D_{\mathcal{A}}(p,q)-1)^{2}+\cdots$$
which is obviously well-defined and is a normalized Dedekind symbol. The corresponding result holds for multiple reciprocity functions.

From this example, by multiplying multiple reciprocity functions above, we have the following trivial result:

\begin{prop}\label{product to reciprocity functions}
Given any indexed set $\mathcal{A}=\{A_{1},A_{2}\cdots\}$ and a family of $R$-valued reciprocity functions $F^{A_{1}},F^{A_{2}},\cdots:U\rightarrow R$. There is an $R$-valued multiple reciprocity function whose $A_{i}$-component is the function $F^{A_{i}}(p,q), \forall A_{i}\in \mathcal{A}$.
\end{prop}

The above examples and proposition inspire us to define:

\begin{Def}\label{Def of indecom MDS}
Given a multiple Dedekind symbol $D_{\mathcal{A}}$. For any two disjoint subsets $\mathcal{A}_{1},\mathcal{A}_{2}\subset\mathcal{A}$, define two multiple Dedekind symbols
$$D_{i}=1+\sum_{w\in \mathcal{A}_{i}}D^{w}w, \ i=1,2,$$
if the equation $D_{\mathcal{A}}=D_{1}D_{2}$ implies $D_{1}$ or $D_{2}$ is trivial, we call $D_{\mathcal{A}}$ indecomposable.
\end{Def}

\subsection{Shuffle property}

In this subsection, we consider the shuffle property of multiple Dedekind symbols and multiple reciprocity functions. For two words
$$u=A_{1}\cdots A_{r}, v=A_{r+1}\cdots A_{r+s}\in \mathcal{A}^{*},$$
denote the shuffle set of $u, v$ by
$$Sh(u,v)=\{w=A_{\sigma^{-1}(1)}\cdots A_{\sigma^{-1}(r+s)}\in \mathcal{A}^{*}\},$$
where $\sigma\in S(r+s)$ is a permutation of order $r+s$ satisfying that
$$\sigma(1)<\cdots<\sigma(r), \sigma(r+1)<\cdots<\sigma(r+s).$$

\begin{Def}
A multiple Dedekind symbol $D_{\mathcal{A}}=1+\sum\limits_{\emptyset\neq w\in \mathcal{A}^{*}}D^{w}w$ is called shuffled if for any $u, v\in A^{*}$, we have
$$D^{u}(p,q)D^{v}(p,q)=\sum_{w\in Sh(u,v)}D^{w}(p,q),\ \forall (p,q)\in U.$$
Replace $D_{\mathcal{A}}$ by $F_{\mathcal{A}}$, and $D^{w}$ by $F^{w}$, we can also define shuffled multiple reciprocity functions.
\end{Def}

Obviously, most of multiple Dedekind symbols are not shuffled. For instance, given two non-trivial shuffled multiple Dedekind symbols, the sum of them minus $1$ is still a multiple Dedekind symbols but not shuffled in general. The simplest shuffled multiple Dedekind symbol is as below:

\begin{exa}
For an indexed set $\mathcal{A}=\{a\}$ and given $\mathbb{C}$-valued Dedekind symbol $D(p,q)$, the multiple Dedekind symbol $D_{\mathcal{A}}(p,q)=e^{D(p,q)a}$ constructed in Section $4.1$ is shuffled. A corresponding result holds for the multiple reciprocity functions.
\end{exa}

If fact, if $R$ is a $k$-algebra, then a multiple Dedekind symbol $D_{\mathcal{A}}$ is shuffled if and only if for any $(p,q)\in U$, the element $D_{\mathcal{A}}(p,q)\in R\langle\langle \mathcal{A}\rangle\rangle^{\times}$ is a group-like element, where we provide $R\langle\langle \mathcal{A}\rangle\rangle$ the natural Hopf algebra structure. Thus one can prove the following Theorem \ref{Shuffled equivalence} and Proposition \ref{product keep shuffled} in a similar but much simpler way. However, in order to avoid the abstract discussion about Hopf algebras, we give a self-contained and more explicit proof here.

The following lemma is straightforward:

\begin{lem}\label{shuffle decom}
Given a multiple Dedekind symbol and its associated multiple reciprocity function
$$D_{\mathcal{A}}(p,q)=1+\sum\limits_{\emptyset\neq w\in \mathcal{A}^{*}}D^{w}(p,q)w,\ F_{\mathcal{A}}(p,q)=1+\sum\limits_{\emptyset\neq w\in \mathcal{A}^{*}}F^{w}(p,q)w.$$
For any $\emptyset\neq u,v\in \mathcal{A}^{*}$, we have
$$\sum_{w\in Sh(u,v)}D^{w}(p,q)=\sum_{u_{1}u_{2}=u,v_{1}v_{2}=v}[\sum_{w_{1}\in Sh(u_{1},v_{1})}F^{w_{1}}(p,q)
\sum_{w_{2}\in Sh(u_{2},v_{2})}D^{w_{2}}(-q,p)].$$
\end{lem}

\begin{Thm}\label{Shuffled equivalence}
A normalized multiple Dedekind symbol $D_{\mathcal{A}}$ is shuffled if and only if its associated multiple reciprocity function $F_{\mathcal{A}}$ is shuffled.
\end{Thm}

\noindent{\bf Proof:}
We only prove the sufficiency, the necessity is similar and easier. Precisely, we need to show that for any $\emptyset\neq u,v\in \mathcal{A}^{*}$,
$$F^{u}(p,q)F^{v}(p,q)=\sum_{w\in Sh(u,v)}F^{w}(p,q)\ \text{implies}\ D^{u}(p,q)D^{v}(p,q)=\sum_{w\in Sh(u,v)}D^{w}(p,q).$$
We prove it by induction on $n=l(u)+l(v)$.

When $n=2$, $u,v\in \mathcal{A}^{*}$ are of length $1$, by the construction of normalized generalized Dedekind symbols from reciprocity functions in Section $3$, we have
\begin{equation*}
\begin{split}
& \ \ \ D^{u}(p,q)D^{v}(p,q)\\
& =\sum\limits_{i=1}^{n}F^{u}(p_{i},q_{i})\sum\limits_{i=1}^{n}F^{v}(p_{i},q_{i})\\
& =\sum\limits_{i=1}^{n}F^{uv}(p_{i},q_{i})+F^{vu}(p_{i},q_{i})+(\sum\limits_{1\leq i<j\leq n}+\sum\limits_{1\leq j<i\leq n})F^{u}(p_{i},q_{i})F^{v}(p_{j},q_{j}),
\end{split}
\end{equation*}
where $(p_{i},q_{i})$ for $i=1,\cdots,n$ are determined by a given continued fraction representation of $(p,q)$. On the other hand, we have
\begin{equation*}
\begin{split}
& \ \ \ D^{uv}(p,q)+D^{vu}(p,q)\\
& =\sum\limits_{i=1}^{n}F^{uv}(p_{i},q_{i})+F^{vu}(p_{i},q_{i})-F^{u}(p_{i},q_{i})D^{v}(p_{i},q_{i})-F^{v}(p_{i},q_{i})D^{u}(p_{i},q_{i}).
\end{split}
\end{equation*}
Notice that the continued fraction representation of $(p_{i},q_{i})$ is determined by $(p_{j},q_{j})$ for $j>i$, by the construction from reciprocity functions to normalized Dedekind symbols, one has
$$D^{u}(p_{i},q_{i})=-\sum_{j>i}F^{u}(p_{j},q_{j}).$$
Thus $D^{u}(p,q)D^{v}(p,q)=D^{uv}(p,q)+D^{vu}(p,q)$.

Now we assume that our statement holds for $l(u)+l(v)\leq n-1$, we consider the case of $n$ below. According to the law $(MDS\ 3)$ and Lemma \ref{shuffle decom}, we have
$$D^{u}(p,q)D^{v}(p,q)=\sum_{u_{1}u_{2}=u,v_{1}v_{2}=v}F^{u_{1}}(p,q)D^{u_{2}}(-q,p)F^{v_{1}}(p,q)D^{v_{2}}(-q,p)$$
and
$$\sum_{w\in Sh(u,v)}D^{w}(p,q)=\sum_{u_{1}u_{2}=u,v_{1}v_{2}=v}[\sum_{w_{1}\in Sh(u_{1},v_{1})}F^{w_{1}}(p,q)\sum_{w_{2}\in Sh(u_{2},v_{2})}D^{w_{2}}(-q,p)].$$
By induction, we have
$$D^{u}(p,q)D^{v}(p,q)-\sum_{w\in Sh(u,v)}D^{w}(p,q)=D^{u}(-q,p)D^{v}(-q,p)-\sum_{w\in Sh(u,v)}D^{w}(-q,p).$$
Denote by $E_{u,v}(p,q)= D^{u}(p,q)D^{v}(p,q)-\sum\limits_{w\in Sh(u,v)}D^{w}(p,q) $ for convenience, then we have
$$E_{u,v}(p,q)=E_{u,v}(-q,p), E_{u,v}(p,p+q)=E_{u,v}(p,q).$$
Thus $E_{u,v}:U\rightarrow R$ is a constant function. Since $D_{\mathcal{A}}(p,q)$ is normalized and furthermore $D^{w}(1,1)=0$ for any $\emptyset\neq w\in \mathcal{A}^{*}$, we have $E_{u,v}(1,1)=0$. Thus $E_{u,v}(p,q)=0$ for any $(p,q)\in W$. It follows that
$$D^{u}(p,q)D^{v}(p,q)=\sum_{w\in Sh(u,v)}D^{w}(p,q).$$
The theorem holds.
$\hfill\Box$\\

\begin{rem}
When $D_{\mathcal{A}}$ is not normalized, we still have the necessity with the same proof. But the converse statement does not hold in general.
\end{rem}

Combining with the following proposition, we finish the proof of Theorem \ref{Intro: shuffle and product} in the introduction.

\begin{prop}\label{product keep shuffled}
Suppose that $D_{\mathcal{A}}$, $E_{\mathcal{A}}$ are two $R$-valued shuffled multiple Dedekind symbols. Then the inverse $D_{\mathcal{A}}^{-1}$ and the product $D_{\mathcal{A}}E_{\mathcal{A}}$ are $R$-valued shuffled multiple Dedekind symbols.
\end{prop}

\noindent{\bf Proof:}
Suppose that
$$D_{\mathcal{A}}(p,q)=1+\sum_{\emptyset\neq w\in \mathcal{A}^{*}}D^{w}(p,q)w,
\ E_{\mathcal{A}}(p,q)=1+\sum_{\emptyset\neq w\in \mathcal{A}^{*}}E^{w}(p,q)w$$
are two shuffled multiple Dedekind symbols. Obviously, $K_{\mathcal{A}}(p,q)=D^{-1}_{\mathcal{A}}(p,q)$ and $H_{\mathcal{A}}(p,q)=D_{\mathcal{A}}(p,q)E_{\mathcal{A}}(p,q)$ are $R$-valued multiple Dedekind symbols.

For the first statement, notice that for any word $u,v\in \mathcal{A}^{*}$ of length $1$, we have $K^{u}(p,q)=-D^{u}(p,q)$ and $K^{uv}(p,q)=-D^{uv}(p,q)+D^{u}(p,q)D^{v}(p,q)$. It is trivial to check that
$$K^{u}(p,q)K^{v}(p,q)=K^{uv}(p,q)+K^{vu}(p,q).$$

Suppose that for any $u,v\in \mathcal{A}^{*}$, when $l(u)+l(v)<n$ the statement holds, we want to prove this in the case that $l(u)+l(v)=n+1$. According to Lemma \ref{shuffle decom}, omit the $(p,q)$ in the formula for simplicity, we have
\begin{equation*}
\begin{split}
& 0=\sum\limits_{u=u_{1}u_{2},v=v_{1}v_{2}}\sum\limits_{w_{1}\in Sh(u_{1},v_{1})}\sum\limits_{w_{2}\in Sh(u_{2},v_{2})}D^{w_{1}}E^{w_{2}}\\
& \ \ \ =\sum\limits_{w\in Sh(u,v)}E^{w}-E^{u}E^{v}+\sum\limits_{u=u_{1}u_{2},v=v_{1}v_{2}}D^{u_{1}}D^{v_{1}}E^{u_{2}}E^{v_{2}}\\
& \ \ \ =\sum\limits_{w\in Sh(u,v)}E^{w}-E^{u}E^{v}+(\sum\limits_{u=u_{1}u_{2}}D^{u_{1}}E^{u_{2}})(\sum\limits_{v=v_{1}v_{2}}D^{v_{1}}E^{v_{2}})\\
& \ \ \ =\sum\limits_{w\in Sh(u,v)}E^{w}-E^{u}E^{v}
\end{split}
\end{equation*}
by induction, thus we proved the first statement.

Next, we consider the second statement. For any words $w_{1}$ and $w_{2}$ in $\mathcal{A}^{*}$, we have
$$H^{w_{i}}(p,q)=\sum_{u_{i}v_{i}=w_{i}}D^{u_{i}}(p,q)E^{v_{i}}(p,q), i\in\{1,2\}.$$
Then according to Lemma \ref{shuffle decom}, we have
\begin{equation*}
\begin{split}
H^{w_{1}}H^{w_{2}}
& =\sum\limits_{u_{1}v_{1}=w_{1}}D^{u_{1}}E^{v_{1}}
\sum\limits_{u_{2}v_{2}=w_{2}}D^{u_{2}}E^{v_{2}}\\
& =\sum\limits_{u_{1}v_{1}=w_{1},u_{2}v_{2}=w_{2}}
\sum\limits_{x\in Sh(u_{1},u_{2})}\sum\limits_{y\in Sh(v_{1},v_{2})}D^{x}E^{y}\\
& =\sum\limits_{w\in Sh(w_{1},w_{2})}
\sum\limits_{xy=w}D^{x}E^{y}=\sum\limits_{w\in Sh(w_{1},w_{2})}H^{w}.
\end{split}
\end{equation*}
Thus the proposition holds.
$\hfill\Box$\\

\section{Relationship to Modular Forms}

In this section we construct and study $\mathbb{C}$-valued shuffled multiple Dedekind symbols and multiple reciprocity functions from modular forms by regularized iterated integrals.

\subsection{Regularized iterated integrals}

In this and the next section, let $\Gamma=SL_{2}(\mathbb{Z})$ and denote by $\mathcal{M}_{2k+2}(\Gamma)$ the $\mathbb{C}$-vector space of modular forms for $\Gamma$ of weight $2k+2$. In particular, it includes the Eisenstein seires of weight $2k+2$.

Recall notations in the introduction, suppose that we are given a set of non-commutative symbols $\mathcal{A}=\{A_{1},A_{2},\cdots\}$ as the indexed set, and that every element $A_{i}$ is associated to a positive even integer $w_{i}$. For any word $B=A_{1}\cdots A_{l}$ of symbols of $\mathcal{A}$, denote by $w(B)=w_{1}+\cdots+w_{l}$. Define $\mathcal{B}(\mathcal{A})\subset\mathcal{A}^{*}$ to be the subset of non-empty words. Denote by
$$\mathcal{M}(\mathcal{A},\Gamma)=\prod_{B\in \mathcal{B}(\mathcal{A})}\mathcal{M}_{w(B)+2}(\Gamma).$$
For any element $\mathbf{h}\in \mathcal{M}(\mathcal{A},\Gamma)$, denote by $\mathbf{h}(B)$ the component of $\mathbf{h}$ at the factor $B\in \mathcal{B}(\mathcal{A})$, it is a modular form of weight $w(B)+2$.

Given an element $\mathbf{h}\in \mathcal{M}(\mathcal{A},\Gamma)$, denote by $\mathfrak{h}$ the upper half complex plane. For any $\tau\in\mathfrak{h}$, consider the differential form
$$\Omega_{\mathbf{h}}(\tau)=\sum_{B\in\mathcal{B}(\mathcal{A})}B\mathbf{h}(B)(X-Y\tau)^{w(B)}d\tau.$$
It is a (formal) differential form whose coefficients are products of symbols $B\in \mathcal{B}(\mathcal{A})$ and homogeneous polynomials of $X$ and $Y$.

There is a natural action of $\Gamma$ on the homogeneous polynomial of $X$ and $Y$. Precisely,
$$\gamma\left(
                   \begin{array}{c}
                     X \\
                     Y \\
                   \end{array}
                 \right)
=\left(
   \begin{array}{c}
     aX+bY \\
     cX+dY \\
   \end{array}
 \right)
,\forall \gamma=\left(
                                                       \begin{array}{cc}
                                                         a & b \\
                                                         c & d \\
                                                       \end{array}
                                                     \right)
\in \Gamma.$$
Denote the action of $\gamma$ on the homogeneous polynomial $P(X,Y)$ by $P(X,Y)|_{\gamma}$. There is also a natural action of $\Gamma$ on the complex upper half plane as
$$\gamma(\tau)=\frac{a\tau+b}{c\tau+d},\ \forall \tau\in\mathfrak{h},\gamma=\begin{pmatrix}
	a & b\\
	c & d
\end{pmatrix}\in \Gamma.$$
Then one find the differential form $\Omega_{\mathbf{h}}(\tau)$ is $\Gamma$-invariant.

It is obvious to see that $d\Omega_{\mathbf{h}}(\tau)+\Omega_{\mathbf{h}}(\tau)\wedge \Omega_{\mathbf{h}}(\tau)=0$. Since the complex upper half plane $\mathfrak{h}$ is simply connected, for any path $\gamma:[0,1]\rightarrow \mathfrak{h}$ the iterated integral
$$I_{\gamma}(\Omega_{\mathbf{h}}(\tau))=1+\int_{\gamma}\Omega_{\mathbf{h}}(\tau)+\int_{\gamma}\Omega_{\mathbf{h}}(\tau)\Omega_{\mathbf{h}}(\tau)+\cdots$$
only depends on the start point $\gamma_{0}=\gamma(0)$ and the end point $\gamma_{1}=\gamma(1)$ of the path $\gamma$, thus we will denote the iterated integral above by $I_{\mathbf{h}}(\gamma_{0},\gamma_{1})$ for convenience.

Now we extend the definition of the iterated integrals $I_{\mathbf{h}}(\tau_{0},\tau_{1})$ to the case that $\tau_{0},\tau_{1}\in\overline{\mathfrak{h}}=\mathfrak{h}\cup\mathbb{Q}\cup\{i\infty\}$. We use the technology of Brown \cite{F.B}. Define
$$\Omega^{\infty}_{\mathbf{h}}(\tau)=\sum_{B\in\mathcal{B}(\mathcal{A})}B\mathbf{h}^{\infty}(B)(X-Y\tau)^{w(B)}d\tau,$$
where $\mathbf{h}^{\infty}(B)$ is the constant term of the Fourier expansion  of the modular form $\mathbf{h}(B)$. We also have $d \Omega_{\mathbf{h}}^{\infty}(\tau) +\Omega_{\mathbf{h}}^{\infty}(\tau)\wedge \Omega_{\mathbf{h}}^{\infty}(\tau)=0$. Similar to the above discussion, we can define the iterated integral
$$I_{\mathbf{h}}^{\infty}(\tau_{0},\tau_{1})=1+\int_{\tau_{0}}^{\tau_{1}} \Omega_{\mathbf{h}}^{\infty}(\tau) +\int_{\tau_{0}}^{\tau_{1}} \Omega_{\mathbf{h}}^{\infty}(\tau) +\cdots.$$
It does not depend on the choice of the path $\gamma$ from $\tau_{0}$ to $\tau_{1}$. Denote by
$$RI_{\mathbf{h}}(\tau,x)=I_{\mathbf{h}}(\tau,x)I_{\mathbf{h}}^{\infty}(x,\tau),\forall \tau,x\in\mathfrak{h}.$$

\begin{lem}[Brown \cite{F.B}]
$RI_{\mathbf{h}}(\tau,x)$ is finite as $x\rightarrow i\infty$ and converges like $\mathcal{O}(e^{2\pi ix})$.
\end{lem}

Now we can define the iterated Eichler integral from $\tau\in\mathfrak{h}$ to the tangential base point $\overrightarrow{s}_{\infty}$ at the cusp point $i\infty$, where $s\in\mathbb{Q}$ is a rational number. By the above lemma, we can define:

\begin{Def}
Denote by
$$I_{\mathbf{h}}(\tau,\overrightarrow{s}_{\infty})=\lim_{\varepsilon\rightarrow i\infty}I_{\mathbf{h}}(\tau,\varepsilon)I_{\mathbf{h}}^{\infty}(\varepsilon,s).$$
We call it the iterated Eichler integral of $\mathbf{h}\in\mathcal{M}(\mathcal{A};\Gamma)$ from $\tau$ to the tangential base point $\overrightarrow{s}_{\infty}$. We also call the coefficients of symbols $B\in\mathcal{B}(\mathcal{A})$ regularized iterated integrals of $\mathbf{h}(A_{1}),\cdots,\mathbf{h}(A_{2})$, where $B=A_{1}\cdots A_{2}$.
\end{Def}

One finds that $I_{\mathbf{h}}(\tau,\overrightarrow{s}_{\infty})$ is a function of $\tau$, whose coefficients are words of $\mathcal{A}^{*}$ timing Laurent series of $X$ and $Y$.

\begin{rem}
The technologies and definitions above are from Brown \cite{F.B}, but our notation of tangential base point is different from Brown's one.
\end{rem}

It is natural to regard $\mathcal{M}_{\mathcal{A}}(\Gamma)=\prod\limits_{A\in\mathcal{A}}\mathcal{M}_{w(A)+2}(\Gamma)$ as a subspace of $\mathcal{M}(\mathcal{A},\Gamma)$, i.e, the element $\mathbf{f}=(f_{1},f_{2},\cdots)\in \mathcal{M}_{\mathcal{A}}(\Gamma)$ can be regarded as the element $\mathbf{h}\in \mathcal{M}(\mathcal{A},\Gamma)$ by
$$\mathbf{h}(A_{i})=f_{i}\ \text{and}\ \mathbf{h}(B)=0,\forall d(B)\geq2.$$
Take the rational point $s=0$ and $\mathbf{h}\in \mathcal{M}_{\mathcal{A}}(\Gamma)$, then the iterated Eichler integral $I_{\mathbf{h}}(\tau,\overrightarrow{0}_{\infty})$ is the one defined by Brown \cite{F.B}, denote it by $I_{\mathbf{h}}(\tau,\infty)$ for convenience.

We can also define iterated Eichler integrals when $\tau\in\mathbb{Q}$. Precisely, for any $\tau_{1}\in\mathfrak{h}$ define
$$I_{\mathbf{h}}(\overrightarrow{s}_{\tau},\tau_{1})=I_{\mathbf{h}}(\overrightarrow{\gamma^{-1}(s)}_{\gamma^{-1}(\tau)},\gamma^{-1}(\tau_{1}))|_{\gamma},$$
where $\gamma\in\Gamma$ is the element satisfying $\gamma(\infty)=\tau$. Furthermore, we can define
$$I_{\mathbf{h}}(\overrightarrow{s}_{\tau_{0}},\overrightarrow{t}_{\tau_{2}})=I_{\mathbf{h}}(\overrightarrow{s}_{\tau_{0}},\tau_{1})I_{\mathbf{h}}(\tau_{1},\overrightarrow{t}_{\tau_{2}}),$$
where $\overrightarrow{s}_{\tau_{0}}$ and $\overrightarrow{t}_{\tau_{2}}$ are tangential base points and $\tau_{1}\in\mathfrak{h}$ is an arbitrary point. It is independent of the choice of $\tau_{1}$ and thus well-defined according to the following proposition, one can find the proof in Proposition $4.7$ of \cite{F.B}.
\begin{prop}\label{II property}
The iterated Eichler integral satisfies the following properties:

$(1)$. (Differential property). $\frac{d}{d\tau}I_{\mathbf{h}}(\tau,\overrightarrow{s}_{\infty})=-\Omega(\tau)I_{\mathbf{h}}(\tau,\overrightarrow{s}_{\infty}),\ \forall\tau\in\mathfrak{h}$.

$(2)$. (Product). $I_{\mathbf{h}}(\tau,\overrightarrow{s}_{\infty})=I_{\mathbf{h}}(\tau,\tau_{1}) I_{\mathbf{h}}(\tau_{1},\overrightarrow{s}_{\infty}),\ \forall\tau_{1}\in\mathfrak{h}$.

$(3)$. (Shuffle). $I_{\mathbf{h}}(\tau,\overrightarrow{s}_{\infty})$ is group-like and invertible.

$(4)$. (Modular). $I_{\mathbf{h}}(\gamma(\tau),\overrightarrow{\gamma(s)}_{\gamma(\infty)})|_{\gamma}=I_{\mathbf{h}}(\tau,\overrightarrow{s}_{\infty}),\ \forall\gamma\in\Gamma$.
\end{prop}

The shuffle property and modular property are particularly important for our later construction.

\subsection{Constructions from modular forms}

With the help of the iterated Eichler integrals, we can construct a family of multiple Dedekind symbols and the associated multiple reciprocity functions explicitly from modular forms for the group $\Gamma$. For $\mathcal{A}$ and $\mathbf{h}$ as above, denote by
$$D_{\mathbf{h}}(p,q)=I_{\mathbf{h}}(\overrightarrow{q/p}_{\infty},\overrightarrow{\infty}_{q/p})(p,q), F_{\mathbf{h}}(p,q)=I_{\mathbf{h}}(\overrightarrow{q/p}_{\infty},\overrightarrow{q/p}_{0})(p,q)$$
and
$$E_{\mathbf{h}}(p,q)=\text{exp}(\sum_{A_{i}\in\mathcal{A}}A_{i}a_{\mathbf{h}(A_{i})}(0)/pq)=1+\sum\frac{a_{\mathbf{h}(A_{i_{1}})}(0)\cdots a_{\mathbf{h}(A_{i_{r}})}(0)}{r!(pq)^{r}}A_{i_{1}}\cdots A_{i_{r}},$$
where the iterated integrals are defined on $(p,q)\in U$ by replacing $(p,q)$ by $(Y,X)$, and $a_{\mathbf{h}(A_{i})}(0)$ is the constant term of the Fourier expression of the modular form $\mathbf{h}(A_{i})$. Then we have:

\begin{Thm}\label{modular forms to MRF}
The function $D_{\mathbf{h}}$ is a $\mathbb{C}$-valued shuffled multiple Dedekind symbol based on $\mathcal{A}$, the functions $E_{\mathbf{h}}$ and $F_{\mathbf{h}}$ are $\mathbb{C}$-valued shuffled multiple reciprocity functions based on $\mathcal{A}$. Furthermore, we have
$$F_{\mathbf{h}}(p,q)=D_{\mathbf{h}}(p,q)E_{\mathbf{h}}(p,q)D_{\mathbf{h}}^{-1}(-q,p).$$
\end{Thm}

\noindent{\bf Proof:}
It is obvious to see that $E_{\mathbf{h}}$ is a $\mathbb{C}$-valued shuffled multiple reciprocity function based on $\mathcal{A}$, and directly form the modular property and product property of iterated Eichler integrals, we have
$$D_{\mathbf{h}}(p,q)E_{\mathbf{h}}(p,q)D_{\mathbf{h}}^{-1}(-q,p)=F_{\mathbf{h}}(p,q),\forall (p,q)\in U.$$
Thus we only need to show that $D_{\mathbf{h}}$ is a $\mathbb{C}$-valued shuffled multiple Dedekind symbol based on $\mathcal{A}$, i.e., we need to show that for any $(p,q)\in U$,
$$D_{\mathbf{h}}(p,-q)=D_{\mathbf{h}}(-p,q),$$
which is obvious, and that $D_{\mathbf{h}}(p,q)=D_{\mathbf{h}}(p,p+q)$.

Take $T=\begin{pmatrix}
	1 & 1\\
	0 & 1
\end{pmatrix}\in\Gamma
$ and consider the modular property of iterated Eichler integrals
$$I_{\mathbf{h}}(T(\tau_{0}),T(\tau_{1}))|_{T}=I_{\mathbf{h}}(\tau_{0},\tau_{1}),$$
we have the following equations:
\begin{equation*}
\begin{split}
D_{\mathbf{h}}(p,p+q)& =I_{\mathbf{h}}(\overrightarrow{(p+q)/p}_{\infty},\overrightarrow{\infty}_{(p+q)/p})(p,p+q)\\
& =I_{\mathbf{h}}(\overrightarrow{q/p}_{\infty},\overrightarrow{\infty}_{q/p})(p,q)\\
& =D_{\mathbf{h}}(p,q).
\end{split}
\end{equation*}
Thus $D_{\mathbf{h}}$ is a $\mathbb{C}$-valued multiple Dedekind symbol. In particular, from the shuffle property of iterated Eichler integrals, $D_{\mathbf{h}}$ is shuffled.
$\hfill\Box$\\

\begin{rem}
When there are Eisenstein series in $\mathbf{h}(B)$, the multiple Dedekind symbols and multiple reciprocity functions constructed here are almost.
\end{rem}

Denote by $M_{\mathbf{h}}$ the multiple reciprocity function associated to $D_{\mathbf{h}}$. If we choose an element $\mathbf{h}\in \mathcal{M}_{\mathcal{A}}(\Gamma)$ such that for any $A_{i}\in\mathcal{A}$, $\mathbf{h}(A_{i})$ is a cusp form. Then the multiple Dedekind symbols $D_{\mathbf{h}}$ and the multiple reciprocity function $F_{\mathbf{h}}$ is just the one constructed by Manin \cite{Y.M3}, and $M_{\mathbf{h}}=F_{\mathbf{h}}$ in this case.

It is obvious to see that the multiple reciprocity functions $F_{\mathbf{h}}$ of modular forms are all in Laurent form, which means for any non-empty word $w$ and $(p,q)\in U$, its $w$ component is in the form $G_w(p,q)/(pq)^{n}$ for a positive integer $n$ and a polynomial $G_w$ not depending on $(p,q)$. Moreover, we have $n\leq l(w)$. However, the statement does not hold for the multiple Dedekind symbol $D_{\mathbf{h}}$. In particular, for any $\mathbf{h}_{1},\mathbf{h}_{2}\in\mathcal{M}(\mathcal{A},\Gamma)$, the product of multiple reciprocity functions $F_{\mathbf{h}_{1}}\bullet F_{\mathbf{h}_{2}}$ is not in Laurent polynomial form in general. We give a conjecture here:

\begin{conj}
Suppose that $\mathbf{h}\in\mathcal{M}_{\mathcal{A}}(\Gamma)$ and that there exists at least two elements $A_{1}, A_{2}\in\mathcal{A}$ such that $\mathbf{h}(A_{i})$ is non-trivial. Then the Dedekind symbol $D_{\mathbf{h}}(p,q)$ is indecomposable.
\end{conj}

The conjecture holds if for any family of non-trivial modular forms $f_{1},\cdots, f_{n}$, $n\geq2$, and $(p,q)\in U$ one has
$$D_{f_{1},\cdots,f_{n}}(p,q)\neq D_{f_{1}}(p,q)D_{f_{2},\cdots,f_{n}}(p,q),$$
where $D_{f_{1},\cdots,f_{n}}(p,q)$ is the regularized iterated integral of $f_{i}(\tau)(p\tau-q)^{2k_{i}-2}d\tau$ for $\overrightarrow{\infty}_{q/p}$ to $\overrightarrow{q/p}_{\infty}$ with $2k_{i}$ the weight of $f_{i}$.

\subsection{Length $2$ components associated to Eisenstein series}

For a multiple Dedekind symbol of modular forms, its length one components associated to Eisenstein series were calculated by Fukuhara in \cite{S.F2} and \cite{S.F5}, and we showed the reciprocity law in Formula (\ref{length 1 reciprocity law}). Now we consider length two components.

First we recall the holomorphic multiple modular values defined by Brown \cite{F.B}. Let $\mathcal{B}$ be a rational basis of $\bigoplus_{k\geq1}\mathcal{M}_{2k+2}(\Gamma)$. We take $\mathcal{A}$ to be the set of symbols $\{A_{f};f\in\mathcal{B}\}$, in where $A_{f}$ corresponds to the weight of the modular form $f$. Now take the element $\mathbf{h}\in \mathcal{M}_{\mathcal{A}}(\Gamma)\subset \mathcal{M}(\mathcal{A},\Gamma)$ such that $\mathbf{h}(A_{f})=(2\pi i)^{2k-1}f$ for $f\in\mathcal{A}$ of weight $2k$. Denote by $V_{2k-2}$ the $\mathbb{C}$-vector space of homogeneous polynomial of $X$, $Y$ of degree $2k-2$, and by
$$M_{2k}^{\vee}=\langle A_{f}; f\in\mathcal{A}, w_{f}=2k\rangle,
M^{\vee}=\bigoplus_{k\geq1} M_{2k}^{\vee}\otimes V_{2k-2}.$$
Then Brown \cite{F.B} proved the following proposition:

\begin{prop}
For any $\gamma\in \Gamma$, there exists a series $\mathcal{C}_{\gamma}\in \mathbb{C}\langle\langle M^{\vee}\rangle\rangle$ such that
$$I_{\mathbf{h}}(\tau,\infty)=I_{\mathbf{h}}(\gamma(\tau),\infty)|_{\gamma}\mathcal{C}_{\gamma}.$$
The series $\mathcal{C}_{\gamma}$ does not depend on the choice of $\tau\in\mathfrak{h}$. Moreover, $\mathcal{C}$ determines a (non-commutative) cocycle in $Z^{1}(\Gamma,\mathbb{C}\langle\langle M^{\vee}\rangle\rangle)$.
\end{prop}

\begin{Def}
Define the ring of (holomorphic) multiple modular values $\mathcal{MMV}^{hol}$ for $\Gamma$ to be the $\mathbb{Q}$-algebra generated by the coefficients of $\mathcal{C}_{\gamma}$ for any $\gamma\in \Gamma$.
\end{Def}

It is well-known that the group $\Gamma$ is generated by two matrices $$S=\left(
                                                         \begin{array}{cc}
                                                           0 & -1 \\
                                                           1 & 0  \\
                                                         \end{array}
                                                       \right)\ \text{and}\
T=\left(
           \begin{array}{cc}
             1 & 1 \\
             0 & 1 \\
           \end{array}
         \right).
$$
Thus the ring $\mathcal{MMV}^{hol}$ only depends on the coefficients of $\mathcal{C}_{S}$ and $\mathcal{C}_{T}$ because of the cocycle properties of $\mathcal{C}$. According to \cite{F.B}, $\mathcal{C}_{T}$ has coefficients in $\mathbb{Q}[2\pi i]$, and since $i\in\mathfrak{h}$ is fixed by $S$, we have the following formula for $\mathcal{C}_S$ by choosing $\tau=i$,
$$\mathcal{C}_S=I_{\mathbf{h}}(i,\infty)|_S^{-1}I_{\mathbf{h}}(i,\infty)=I_{\mathbf{h}}(0,\infty).$$
When replacing $(X,Y)$ by $(q,p)$, we find that $\mathcal{C}_{S}(p,q)$ and $F_{\mathbf{h}}(p,q)$ for the given $\mathbf{h}$ as above only differ from a Laurent power series of $p$ and $q$ because
$$F_{\mathbf{h}}(p,q)=I_{\mathbf{h}}^{\infty}(q/p,0)\mathcal{C}_{S}(p,q)^{-1}I_{\mathbf{h}}^{\infty}(q/p,0)|_{S}^{-1}.$$
Thus length two components of multiple reciprocity functions of Eisenstein series have been determined by Brown \cite{F.B}, in where one finds L-valued of cusp forms outside the critical line.

Next, we consider length two components of the multiple Dedekind symbol of Eisenstein series. Let
$$E_{2k}(\tau)=-\frac{B_{2k}}{4k}+\sum_{n\geq1}\sigma_{2k-1}(n)q^{n}$$
be the Hecke normalized Eisenstein series of weight $2k$. For integers $a,b\geq 2$, denote by $D_{E_{2b}E_{2a}}(p,q)$ the component associated to $E_{2b}(\tau)$ and $E_{2a}(\tau)$ (i.e., the regularized iterated integral of $E_{2b}(\tau)$ and $E_{2a}(\tau)$ as in Section $5.2$). Denote by $Q=e^{2\pi iq/p}$. Similar to the argument in \cite{S.F2}, define
\begin{equation*}
\begin{split}
& S_{1}(r,s)\\
=& \frac{(i)^{r+s}p^{r+s-2}}{(2\pi)^{r+s}}\sum\limits_{n_{1},n_{2}\geq1}\sigma_{2a-1}(n_{1})\sigma_{2b-1}(n_{2})Q^{n_{1}+n_{2}}\int_{0<t_{1}<t_{2}<\infty}e^{-(n_{1}t_{1}+n_{2}t_{2})}t_{1}^{r-1}t_{2}^{s-1}dt_{2}dt_{1}\\
=& \frac{(i)^{r+s}p^{r+s-2}}{(2\pi)^{r+s}}\sum\limits_{n_{1},n_{2}\geq1}\frac{\sigma_{2a-1}(n_{1})\sigma_{2b-1}(n_{2})Q^{n_{1}+n_{2}}}{n_{1}^{r}n_{2}^{s}}\int_{0<t_{1}<\frac{n_{1}}{n_{2}}t_{2}<\infty}e^{-t_{1}-t_{2}}t_{1}^{r-1}t_{2}^{s-1}dt_{2}dt_{1}.
\end{split}
\end{equation*}
For $\lambda>0$, denote by $\Gamma_{\lambda}(r,s)$ as
\begin{equation*}
\begin{split}
\int_{0<t_{1}<\lambda t_{2}<\infty}e^{-t_{1}-t_{2}}t_{1}^{r-1}t_{2}^{s-1}dt_{2}dt_{1}& =\int_{0}^{\infty}(\lambda t_{2}^{r})\Gamma(r)e^{-\lambda t_{2}}\sum_{k\geq0}\frac{(\lambda t_{2})^{k}}{\Gamma(r+k+1)}t_{2}^{s-1}e^{-t_{2}}dt_{2}\\
& =\sum_{k\geq0}\frac{\lambda^{r+k}\Gamma(r)}{\Gamma(r+k+1)}\int_{0}^{\infty}e^{-(\lambda+1)t_{2}}t_{2}^{r+s+k-1}dt_{2}\\
& =\frac{\lambda^{r}\Gamma(r)}{(\lambda+1)^{r+s}}\sum_{k\geq0}\frac{\lambda^{k}}{(\lambda+1)^{k}}\frac{\Gamma(r+k+s)}{\Gamma(r+k+1)}.
\end{split}
\end{equation*}
Then $\Gamma_{\lambda}(r,s)$ is well-defined when $Re(r)$ and $Re(s)$ are positive, and has a meromorphic extension to $r,s\in\mathbb{C}$. Thus we have
$$S_{1}(r,s)=\frac{(i)^{r+s}p^{r+s-2}\Gamma(r)}{(2\pi)^{r+s}}\sum\limits_{n_{1},n_{2}\geq1}\frac{\sigma_{2a-1}(n_{1})\sigma_{2b-1}(n_{2})Q^{n_{1}+n_{2}}}{(n_{1}+n_{2})^{r+s}}\sum_{k\geq0}\frac{n_{1}^{k}}{(n_{1}+n_{2})^{k}}\frac{\Gamma(r+k+s)}{\Gamma(r+k+1)}$$
for $r, s\gg 0$. According to \cite{K.M}, $S_{1}(r,s)$ has a meromorphic extension to $r,s\in \mathbb{C}$. Then by our calculation, $S_{1}(r,s)$ is well-defined at $(2a-1,2b-1)$, we have
\begin{equation}\label{DDS of ES}
D_{E_{2b}E_{2a}}(p,q)=S_{1}(2a-1,2b-1)+S_{2}-S_{3}
\end{equation}
with
\begin{equation*}
\begin{split}
& S_{2}=\frac{(-1)^{a+b}B_{2a}}{4a(2a-1)}p^{2b-3}\frac{\Gamma(2a+2b-2)}{(2\pi)^{2a+2b-2}}\sum\limits_{l=1}^{p}\zeta(2a-1,\frac{l}{p})Li_{2a+2b-2}(Q^{l}),\\
& S_{3}=\frac{(-1)^{a+b}B_{2b}}{4b(2b-1)}p^{2a-3}\frac{\Gamma(2a+2b-2)}{(2\pi)^{2a+2b-2}}\sum\limits_{l=1}^{p}\zeta(2b-1,\frac{l}{p})Li_{2a+2b-2}(Q^{l}),
\end{split}
\end{equation*}
where $\zeta(z,t)$ for $t\in \mathbb{Q}$ is the Hurwitz zeta function, $Li_{2a+2b-2}$ is the $(2a+2b-2)$-th polylogarithm.

\begin{rem}
It seems hard to give an explicit expansion $D_{E_{2a}E_{2b}}(p,q)$. However, by the law $(MDS\ 3)$, we have a reciprocity law of it with the other side connected to double modular values of Eisenstein series.
\end{rem}

\subsection{Decomposition}

In the final part of this section, we give a decomposition of the reciprocity function $F_{\mathbf{h}}$, which help us understand its real or imaginary part. For later convenience of use, we discuss $M_{\mathbf{h}}$ rather than $F_{\mathbf{h}}$. For any finite set $\mathcal{A}$, denote by $\mathcal{T}$ the set of multiple reciprocity functions
$$\{e^{c(D_{f_{1}\cdots f_{n}}(p,q)-D_{f_{1}\cdots f_{n}}(-q,p))w}; c\in\mathbb{C}, w\in\mathcal{B}(\mathcal{A})\},$$
where $f_{i}\in\{\mathbf{h}(B),B\in\mathcal{B}(\mathcal{A})\}$ and $D_{f_{1}\cdots f_{n}}(p,q)$ is the regularized iterated integral of $f_{i}(\tau)(q-p\tau)^{w_{f_{i}}-2}d\tau$ from the tangential base point $\overrightarrow{q/p}_{\infty}$ to $\overrightarrow{\infty}_{q/p}$, it is a $\mathbb{C}$-valued Dedekind symbol by Theorem \ref{modular forms to MRF}.

\begin{Thm}\label{prod decom of MRF}
Suppose that $\mathcal{A}$ is finite. For any $\mathbf{h}\in\mathcal{M}(\mathcal{A},\Gamma)$, there exists a family of elements $T_{i}\in\mathcal{T}$ for $i=0,1,\cdots$ satisfying that

$(1)$. $s_{i}\leq s_{i+1}$, where $s_{i}$ is the maximal positive integer such that the length $k_{i}$ components of $T_{i}$ are zero for $0<k_{i}<s_{i}$.

$(2)$. $M_{\mathbf{h}}(p,q)=(T_{1}\bullet T_{2}\bullet \cdots)(p,q)$ (the product may be infinite).\\
\end{Thm}

\noindent{\bf Proof:}
We give an equivalent proof, i.e., any multiple Dedekind symbol $D_{\mathbf{h}}(p,q)$ can be expressed as a product of elements in the set
$$\mathcal{D}'=\{e^{cD_{f_{1}\cdots f_{n}}(p,q)w},c\in\mathbb{C}, w\in\mathcal{B}(\mathcal{A})\}$$
with $f_{i}\in\{\mathbf{h}(B),B\in\mathcal{B}(\mathcal{A})\}$. We prove this by induction, and use the fact that $D_{\mathbf{h}}(p,q)$ and $e^{cD_{f_{1}\cdots f_{n}}(p,q)w}$ are shuffled.

Suppose that $\mathcal{A}=\{A_{1},\cdots,A_{m}\}$ and that $\mathbf{h}(B)=f_{B}$ for $B\in \mathcal{B}(\mathcal{A})$. Take
$$X_{B}^{(1)}(p,q)=e^{D_{f(B)}(p,q)B}, d(B)=1, B\in\mathcal{B}(\mathcal{A}).$$
Then $D_{\mathbf{h}}(p,q)$ and $\prod\limits_{i} X_{i}^{(1)}(p,q)$ have the same length $1$ components. Now Suppose that we have determine a family of elements $X_{B}^{(k)}$ for $k<n$, $B\in\mathcal{B}(\mathcal{A})$ and $d(B)=k$ satisfying $(1)$ and $D_{\mathbf{h}}(p,q)$ has the same length $<n$ components as
$$D_{n}(p,q)=\prod_{k=1}^{n-1}\prod_{B\in \mathcal{B}(\mathcal{A}), d(B)=k}X_{B}^{(k)}(p,q),$$
where $X_{B}^{(k)}(p,q)$ is a product of elements in $\mathcal{D}'$. Notice that this product depends on the order since the product of multiple Dedekind symbols are not commutative. Any $w$ component of length $n$ of $D_{n}(p,q)$ is a linear combination of $\mathbb{C}$-valued Dedekind symbols in the shape
$$\sum_{i}c'_{i}D_{f_{B_{1}}\cdots f_{B_{t}}}(p,q), B_{i}\in\mathcal{B}(\mathcal{A}), d(B_{i})<n, d(B_{1})+\cdots d(B_{t})=n.$$
Thus by the shuffle relation, it can be written as
$$\sum_{j}c_{j}D_{f_{A_{j_{1}}}\cdots f_{A_{j_{n}}}}(p,q), A_{j_{i}}\in \mathcal{A}.$$
Construct $X_{w}^{(n)}(p,q)$ as
$$X_{w}^{(n)}(p,q)= e^{D_{f_{w}}(p,q)w}\prod_{i}e^{-c_{i}D_{f_{B_{1}}\cdots f_{B_{t}}}(p,q)w}.$$
Then it is obvious that $D_{\mathbf{h}}(p,q)$ has the same length $\leq n$ components as
$$D_{n+1}(p,q)=\prod_{k=1}^{n}\prod_{B\in \mathcal{B}(\mathcal{A}), d(B)=k}X_{B}^{(k)}(p,q),$$
and the decomposition satisfies $(1)$, thus well-defined. Then the theorem follows by induction.
$\hfill\Box$\\

Actually, we have proved the corresponding result for any shuffled multiple reciprocity function. In general, the product is not finite unless the multiple reciprocity function is decomposable, and if we give an order of words of $\mathcal{A}$ (such like the lexicographical order), the decomposition is unique.

Finally we give a question. By Manin \cite{Y.M3}, one could naturally regard a $G_{0}$-valued reciprocity function as a (non-commutative) Dedekind cocycle in $Z_{D}^{1}(\Gamma, G)$, where $G$ is the group of functions defined on $U$ taking values in $G_{0}$.

Take $G_{0}=\mathbb{C}\langle\langle \mathcal{A}\rangle\rangle^{\times}$ with $\mathcal{A}$ as above, for any $c\in H^{1}(\Gamma,G)$, regard it as a multiple reciprocity function. Suppose that the components of $c(p,q)$ are in homogeneous (Laurent) polynomial form, and that $c(p,q)$ is shuffled. The question is that how far does it differ from the cocycle $F_{\mathbf{h}}$ for some $\mathbf{h}\in\mathcal{M}(\mathcal{A},\Gamma)$ as above?

\section{Multiple Dedekind Symbols for Congruence Subgroups}

In the final section, we show how to define multiple Dedekind symbols and multiple reciprocity functions for congruence group $\Gamma\subset SL_{2}(\mathbb{Z})$ of finite index, and state the analogue results.

Keep the assumption in Section $5$, the group $SL_{2}(\mathbb{Z})$ acts on the set $U$ as
$$\gamma:(p,q)\mapsto (cp+dq, ap+bq)$$
with $\gamma=\begin{pmatrix}
a & b \\
c & d
\end{pmatrix}\in SL_{2}(\mathbb{Z})
$. It induces an action of $SL_{2}(\mathbb{Z})$ on Dedekind symbols and reciprocity functions as
$$D|_{\gamma}(p,q)=D(cq+dp,aq+bp), F|_{\gamma}(p,q)=F(cq+dp,aq+bp).$$
Recall that
$$S=\begin{pmatrix}
0 & -1 \\
1 & 0
\end{pmatrix}, T=\begin{pmatrix}
1 & 1 \\
0 & 1
\end{pmatrix},\ V=TS=\begin{pmatrix}
1 & -1 \\
1 & 0
\end{pmatrix}.
$$
Then the properties of Dedekind symbols are
$$D=D|_{-I}, D=D|_{T}$$
and the properties of reciprocity functions are
$$F=F|_{-I},F|_{S}\times F=1,\ F|_{V^{2}}\times F|_{V}\times F=1,$$
where $\times$ is the product in $G$. This inspires that one can define $\mathbb{C}$-valued multiple Dedekind symbols and multiple reciprocity functions for any congruence subgroup $\Gamma\subset SL_{2}(\mathbb{Z})$ of finite index. We will omit the proofs in this section since they are the natural generalization of ones in above sections.

\subsection{Statement}

From now on, let $\Gamma$ be a subgroup of $SL_{2}(\mathbb{Z})$ of finite index. Denote by $G$ the space of $\mathbb{C}$-valued functions in homogeneous (Laurent) polynomial form of $(p,q)\in U$, for any homogenous map
$$f:\Gamma\backslash SL_{2}(\mathbb{Z})\rightarrow G,$$
of degree $2k$, define an action of $SL_{2}(\mathbb{Z})$ on it to be
$$f|_{g}(\gamma)(p,q)=f(\gamma g^{-1})(cp+dq,ap+bq),$$
where $\gamma\in \Gamma\backslash SL_{2}(\mathbb{Z})$ and $g=\begin{pmatrix}
a & b \\
c & d
\end{pmatrix}
\in SL_{2}(\mathbb{Z})$.

Although it is hard to define Dedekind symbols and reciprocity functions for general non-commutative group, we can still define $\mathbb{C}$-valued multiple Dedekind symbols and multiple reciprocity functions via $(MDS\ 1,2)$ and $(MRF\ 1,2,3)$.

\begin{Def}
A $\mathbb{C}$-valued multiple Dedekind symbol based on the set $\mathcal{A}$ for $\Gamma$ is a family of functions
$$D_{\mathcal{A}}=\{D^{w}(\gamma)\in G,w\in \mathcal{A}^{*},\gamma\in \Gamma\backslash SL_{2}(\mathbb{Z})\}$$
satisfying that for any $\gamma\in \Gamma\backslash SL_{2}(\mathbb{Z}), w\in \mathcal{A}^{*}$, $D^{w}(\gamma)$ have the same degree and:

$(MDS\ 1)$. $D^{w}(\gamma)=D^{w}|_{-I}(\gamma)$.

$(MDS\ 2)$. $D^{w}(\gamma)=D^{w}|_{T}(\gamma)$.\\
We also call $D^{w}(\gamma)$ the $w$-component of $D_{A}$ and $l(w)$ its length.
\end{Def}

\begin{Def}
A $\mathbb{C}$-valued multiple reciprocity function based on the set $\mathcal{A}$ is a family of functions
$$F_{\mathcal{A}}=\{F^{w}(\gamma)\in G,w\in \mathcal{A}^{*},\gamma\in \Gamma\backslash SL_{2}(\mathbb{Z})\}$$
satisfying that for any $\gamma\in \Gamma\backslash SL_{2}(\mathbb{Z}), w\in\mathcal{A}^{*}$, $F^{w}(\gamma)$ have the same degree and:

$(MRF\ 1)$. $F^{w}(\gamma)=F^{w}|_{-I}(\gamma)$.

$(MRF\ 2)$. $\sum\limits_{uv=w}F^{u}(\gamma)F^{v}|_{S}(\gamma)=0.$

$(MRF\ 3)$. $\sum\limits_{uv=w}F^{u}|_{V^{2}}(\gamma)F^{v}|_{V}(\gamma)=F^{w}|_{S}(\gamma)$.\\
We also call $F^{w}(\gamma)$ the $w$-component of $F_{\mathcal{A}}$ and $l(w)$ its length.
\end{Def}

\begin{rem}
One can define $\mathbb{C}$-valued Dedekind symbols and reciprocity functions directly in the similar argument.
\end{rem}

The following theorem still hold:

\begin{Thm}\label{DS bijec RF for Gamma}
For any integer $k$, there is a bijection between the set of $\mathbb{C}$-valued normalized (almost) multiple Dedekind symbols and the set of (almost) multiple reciprocity functions of degree $2k$ for $\Gamma$.
\end{Thm}

Also, it is natural to define the product and shuffle property of multiple Dedekind symbols and multiple reciprocity functions by fixing $\gamma\in\Gamma\backslash SL_{2}(\mathbb{Z})$, and the following result holds:

\begin{Thm}
Given a congruence subgroup $\Gamma\subset SL_{2}(\mathbb{Z})$ of finite index and a based set $\mathcal{A}$. The $\mathbb{C}$-valued normalized multiple Dedekind symbols $D_{\mathcal{A}}$ for $\Gamma$ is shuffled if and only if its associated multiple Reciprocity function $F_{\mathcal{A}}$ is shuffled. The product of two $\mathbb{C}$-valued shuffled multiple Dedekind symbols for $\Gamma$ is still shuffled.
\end{Thm}

We can construct $\mathbb{C}$-valued multiple Dedekind symbols and multiple reciprocity functions by modular forms for $\Gamma$. The construction is almost the same, we only need to change the space $\mathcal{M}_{2k}(SL_{2}(\mathbb{Z}))$ to $\mathcal{M}_{2k}(\Gamma)$, and the differential form we need here is
$$\Omega_{\mathbf{h}}(\tau)(\gamma)=\sum_{B\in \mathcal{B}(\mathcal{A})}B\mathbf{h}(B)|_{-(w(B)+2)}\gamma(Y-X\tau)^{w(B)}d\tau,\ \gamma\in \Gamma\backslash SL_{2}(\mathbb{Z}),$$
where $f|_{-2k}\gamma=(c\tau+d)^{-2k}f(\gamma(\tau))$ is the slash operation for any modular form $f$ of weight $2k$. With notations in Section $5$, we have:

\begin{Thm}
The family of functions $\{D_{\mathbf{h}}(\gamma);\gamma\in\Gamma\backslash SL_{2}(\mathbb{Z})\}$ is a $\mathbb{C}$-valued shuffled multiple Dedekind symbol based on $\mathcal{A}$ for $\Gamma$, and the family of functions $\{F_{\mathbf{h}}(\gamma); \gamma\in\Gamma\backslash SL_{2}(\mathbb{Z})\}$ is a $\mathbb{C}$-valued shuffled multiple reciprocity functions based on $\mathcal{A}$ for $\Gamma$.
\end{Thm}

\subsection{The case of $\Gamma_{0}(2)$}

We consider the case of $\Gamma_{0}(2)$ and non-multiple functions as the example. Since
$$\overline{I}=\overline{T}, \overline{V}=\overline{S}, \overline{V^{2}}=\overline{SV}.$$
The axioms for Dedekind symbols are
$$D(\overline{I})(p,q)=D(\overline{I})(-p,-q), D(\overline{V})(p,q)=D(\overline{V})(-p,-q), D(\overline{V^{2}})(p,q)=D(\overline{V^{2}})(-p,-q),$$
$$D(\overline{I})(p,q)=D(\overline{I})(p,p+q), D(\overline{V})(p,q)=D(\overline{V^{2}})(p,p+q), D(\overline{V}^{2})(p,q)=D(\overline{V})(p,p+q).$$
The axioms for reciprocity functions are
$$F(\overline{I})(p,q)=F(\overline{I})(-p,-q), F(\overline{V})(p,q)=F(\overline{V})(-p,-q), F(\overline{V}^{2})(p,q)=F(\overline{V}^{2})(-p,-q),$$
$$F(\overline{I})(p,q)+F(\overline{V})(-q,p)=0, F(\overline{V^{2}})(p,q)+F(\overline{V^{2}})(-q,p)=0$$
$$F(\overline{V^{2}})(q,q-p)+F(\overline{V})(q-p,-p)=F(\overline{S})(-q,p).$$
In fact, the reciprocity function $F$ for $\Gamma_{0}(2)$ are totally determined by $F(\overline{I})$, and the reciprocity law is given by
$$F(\overline{I})(p,q)=D(\overline{I})(p,q)-D(\overline{V})(-q,p),$$
$$F(\overline{V})(p,q)=D(\overline{V})(p,q)-D(\overline{I})(-q,p),$$
$$F(\overline{V^{2}})(p,q)=D(\overline{V^{2}})(p,q)-D(\overline{V^{2}})(-q,p).$$

\begin{rem}
We can give another description for the reciprocity functions for $\Gamma_{0}(2)$ by the Shapiro isomorphism and the Eichler-Shimura isomorphism. One can find more details in \cite{Y.C} and \cite{V.P}. 
\end{rem}

Finally, we calculate the functions $D(\overline{I})$ and $F(\overline{I})$ of Eisenstein series for $\Gamma_{0}(2)$. For $k\geq2$, the space of Eisenstein series of weight $2k$ for $\Gamma_{0}(2)$ has dimension $2$, generated by
$$E_{2k,1}(\tau)=-\frac{B_{2k}}{4k}+\sum_{n\geq1}\sigma_{2k-1}(n)q^{n}\ \text{and}\ E_{2k,2}(\tau)=-\frac{B_{2k}}{4k}+\sum_{n\geq}\sigma_{2k-1}(n)q^{2n}$$
with $q=e^{2\pi i\tau}$. The first one is just the Hecke normalized Eisenstein series for $SL_{2}(\mathbb{Z})$, and we only need to consider the second one. Similar to the argument of Fukuhara \cite{S.F2}, we have
\begin{equation*}
\begin{split}
& \ \ \ D_{E_{2k,2}}(p,q)\\
& =i^{2k-1}p^{2-1}\int_{0}^{\infty}[\sum_{n\geq1}\sum_{m\geq1}n^{2k-1}e^{4\pi imn(it+q/p)}]t^{s-1}dt\\
& =i^{2k-1}p^{s-1}\frac{\Gamma(s)}{(4\pi)^{s}}\sum_{n\geq1}\sum_{m\geq1}\frac{n^{2k-1-s}}{m^{s}}e^{4\pi imnq/p}\\
& =i^{2k-1}p^{2k-2}\frac{\Gamma(s)}{(2\pi)^{s}}\sum_{m=1}^{\infty}\sum_{l=1}^{p}\zeta(s+1-2k,\frac{l}{p})\frac{1}{m^{s}}e^{4\pi imlq/p},\ \text{take}\ s=2k-1\\
& =i^{2k-1}p^{2k-2}\frac{\Gamma(2k-1)}{(4\pi)^{2k-1}}\sum_{l=1}^{p}\zeta(0,l/p)\sum_{m=1}\frac{1}{m^{2k-1}}e^{4\pi imlq/p}.
\end{split}
\end{equation*}
Since $\zeta(0,l/p)=-(\frac{l}{p}-\frac{1}{2})$ and
$$\delta_{2k-1}(p,q)=\sum_{l=1}^{p}\sum_{m=1}^{\infty}\frac{1}{m^{2k-1}}e^{4\pi imlq/p}=\begin{cases}
2^{2k-1}\zeta(2k-1)/p^{2k-2} & p\ even,\\
\zeta(2k-1)/p^{2k-2} & p\ odd,
\end{cases}
$$
we have
$$D_{E_{2k,2}}(\overline{I})(p,q)=\frac{p^{2k-2}(2k-2)!}{(4\pi i)^{2k-1}}[\sum_{l=1}^{p-1}\frac{l}{p}Li_{2k-1}(e^{4\pi ilq/p})+\zeta(2k-1)-\frac{1}{2}\delta_{2k-1}(p,q)].$$
As for the reciprocity function, notice that the L-values of $E_{2k,2}(\tau)$ is
$$L(E_{2k,2},r)=\sum_{n\geq1}\frac{\sigma_{2k-1}(n)}{(2n)^{r}}=2^{-r}L(E_{2k,1},r)$$
Thus
$$F_{E_{2k,2}}(\overline{I})(p,q)=\sum_{r=0}^{2k-2}i^{1-r}{2k-2 \choose r}2^{-r-1}L(E_{2k,1},r+1)p^{r}q^{2k-2-r}$$
$$-\frac{B_{2k}}{4k(2k-1)}(q^{2k-1}p^{-1}+p^{2k-1}q^{-1})-\frac{B_{2k}}{4k}p^{-1}q^{-1}.$$
Similar calculations holds for $F(\gamma)$ for any $\gamma\in \Gamma_{0}(2)\backslash SL_{2}(\mathbb{Z})$.

\end{document}